\documentclass[12pt]{article}

\usepackage{amstext    }
\usepackage{amsthm    }
\usepackage{a4}
\usepackage[mathscr]{eucal}
\usepackage{mathrsfs}

\usepackage{amsmath}
\usepackage{amssymb}
\usepackage{amscd}
\newtheorem{theorem}{Theorem}[section]
\newtheorem{definition}[theorem]{Definition}
\newtheorem{proposition}[theorem]{Proposition}
\newtheorem{corollary}[theorem]{Corollary}
\newtheorem{lemma}[theorem]{Lemma}
\newtheorem{remark}[theorem]{Remark}
\newtheorem{remarks}[theorem]{Remarks}
\newtheorem{example}[theorem]{Example}

\newcommand{\cali}[1]{\mathscr{#1}}

\newcommand{\supp}{{\rm supp}}

\newcommand{\loc}{{loc}}
\newcommand{\ddc}{{dd^c}}

\newcommand{\dbar}{{\overline\partial}}
\newcommand{\ddbar}{{\partial\overline\partial}}

\newcommand{\ind}{{\bf 1}}
\newcommand{\id}{{\rm id}}

\newcommand{\Tan}{{\rm Tan}}

\renewcommand{\Re}{{\rm Re}}
\renewcommand{\Im}{{\rm Im}}

\newcommand{\Cc}{\cali{C}}

\newcommand{\Ec}{\cali{E}}

\newcommand{\C}{\mathbb{C}}
\newcommand{\D}{\mathbb{D}}

\newcommand{\R}{\mathbb{R}}
\newcommand{\T}{\mathbb{T}}
\newcommand{\B}{\mathbb{B}}

\renewcommand{\S}{\mathbb{S}}
\renewcommand{\P}{\mathbb{P}}
\newcommand{\E}{\mathbb{E}}

\title{Density of positive closed currents, a theory of non-generic intersections}

\author{Tien-Cuong Dinh and Nessim Sibony}

\begin{document}

\maketitle

\begin{abstract} 
We introduce a notion of density which extends both the notion of Lelong number  and the theory of intersection for positive closed currents on K\"ahler manifolds. 
For arbitrary finite family of positive closed currents on a compact K\"ahler manifold we construct cohomology classes which represent  their intersection even when a phenomenon of  excess dimension occurs. An example is the case of two algebraic varieties whose intersection has dimension larger than the expected number. The theory allows to solve problems in complex dynamics. Basic calculus on the density of currents is established. 
\end{abstract}

\noindent
{\bf Classification AMS 2010:} 32U, 37F, 32H50


\section{Introduction} \label{intro} 

Positive closed currents are a fundamental tool in complex analysis, algebraic geometry, differential geometry, dynamics and theory of foliations. Bi-degree $(1,1)$ currents and their intersections were intensively studied and had many applications. 
The key point here is that positive closed $(1,1)$-currents can be locally written as the $\ddc$ of plurisubharmonic (psh) functions which are unique up to pluriharmonic functions. Therefore, the study of positive closed $(1,1)$-currents can be systematically reduced to the study of psh functions.
 
The case of arbitrary bi-degree currents is still far from being well-understood. Their local potentials may differ by singular currents and a good choices of potential depend on the choice of coordinates. 
We refer the reader to the works by Demailly, Forn\ae ss, Lelong, M\'eo, Siu, Skoda, Vigny, the authors \cite{Demailly1, DinhSibony3, DinhSibony10,DinhSibony11,FornaessSibony, Lelong,Meo,Sibony2, Siu,Skoda,Vigny2} and the references therein for results and applications  in this direction. The difficulty in the study of higher bi-degree currents is comparable with the study of higher co-dimension cycles in algebraic geometry. The case of bi-degree $(1,1)$ corresponds to the case of hypersurfaces of algebraic varieties.

In \cite[p.16]{Demailly2}, Demailly posed the problem to develop a theory of intersection for positive closed currents of higher bi-degree.
A partial answer to this question was given by the authors  in \cite{DinhSibony10, DinhSibony11} with applications in dynamics, see also Ahn \cite{Ahn} and de Th\'elin-Vigny \cite{deThelinVigny}. In this paper, we develop a new idea which leads to an intrinsic notion of intersection in a quite general setting.
Moreover, the theory of density that we present here covers the notion of Lelong number and the notion of tangent cones of currents at a point. It also permits to measure the size of the intersection for arbitrary currents even when there is an excess of the  intersection dimension. The last phenomenon cannot be studied using the classical intersection theory for currents. See Fulton \cite{Fulton} for the algebraic setting.
 
Let us now, for simplicity, discuss first the case of two positive closed currents $T_1$ and $T_2$ of bidegrees $(p_1,p_1)$ and $(p_2,p_2)$ on a projective manifold $X$ of dimension $k$. 
Consider the particular case where $T_1$ and $T_2$ are given by integration on submanifolds $V_1$ and $V_2$ such that $\dim V_1+\dim V_2<k$. For generic such submanifolds, we have $V_1\cap V_2=\varnothing$. However, in general this intersection may be non-empty and the classical  intersection theory of currents does not give a meaning to this intersection for bi-degree reason.

On the other hand, when $V_2$ is a point, denoted by $a$, there is a notion of multiplicity of $V_1$ at $a$. More generally, if $T_1$ is a general positive closed current there is a notion of Lelong number $\nu(T_1,a)$ of $T_1$ at $a$ which represents the density of $T_1$ at $a$. We first recall this notion and then extend it to the general case. For more detailed expositions on Lelong numbers, see Demailly \cite{Demailly3}, Lelong \cite{Lelong} and Siu \cite{Siu}. 

Choose a local holomorphic coordinate system $x$ near $a$ such that $a=0$ in these coordinates. The Lelong number of $T_1$ at $a$ is the limit of the normalized mass of $T_1$ on the ball $\B(0,r)$ of center 0 and radius $r$ when $r$ tends to 0. More precisely, we have
$$\nu(T_1,a):=\lim_{r\to 0} {\|T_1\|_{\B(0,r)}\over (2\pi)^{k-p_1}r^{2k-2p_1}}\cdot $$
Note that $(2\pi)^{k-p_1}r^{2k-2p_1}$ is the mass on $\B(0,r)$ of the $(p_1,p_1)$-current of integration on a linear subspace of dimension $k-p_1$ through 0. 
Lelong proved that this limit always exists \cite{Lelong}. 
Thie showed that when $T_1$ is given by an analytic set this number is the multiplicity of $V_1$ at $a$. 
Siu proved that the limit does not depend on the choice of coordinates and that the function $a\mapsto \nu(T_1,a)$ is upper semi-continuous for the Zariski topology \cite{Siu}. 

Let $\sigma:\widehat X\to X$ be the blow-up of $X$ at $a$. The pull-back of $T_1$ to $\widehat X\setminus \sigma^{-1}(a)$ is a positive closed current that can be extended by 0 through the exceptional hypersurface $\sigma^{-1}(a)\simeq \P^{k-1}$. We call it the {\it strict transform} of $T_1$ and denote it by $\sigma^\diamond(T_1)$. In general the class of this current in the de Rham cohomology $H^*(\widehat X,\C)$ is not equal to the pull-back by $\sigma$ of the class of $T_1$ in $H^*(X,\C)$. According to Siu's results \cite{Siu}, the missing class can be represented by $\nu(T_1,a)$ times the class of a linear subspace in $\sigma^{-1}(a)$. 

We can also consider the Lelong number from another geometric point of view related to Harvey's exposition \cite{Harvey}. Let $A_\lambda:\C^k\to\C^k$ be defined by $A_\lambda(x):=\lambda x$ with $\lambda\in\C^*$. When $\lambda$ goes to infinity, the domain of definition of the current $T_{1,\lambda}:=(A_\lambda)_*(T_1)$ converges to $\C^k$. 
This family of currents is relatively compact and any limit current, for $\lambda\to\infty$, is invariant under the action of $\C^*$, i.e. invariant under $(A_\lambda)_*$. If $S$ is a limit current, we can extend it to $\P^k$ with zero mass on the hyperplane at infinity. Thus, there is a positive closed current $S_\infty$ on $\P^{k-1}$ such that $S=\pi_\infty^*(S_\infty)$. Here we identify the hyperplane at infinity with $\P^{k-1}$ and we denote by $\pi_\infty:\P^k\setminus \{0\}\to\P^{k-1}$ the canonical central projection (we do not consider the case where $T_1$ is a measure, i.e. $p_1=k$). The class of $S_\infty$ (resp. of $S$) in the de Rham cohomology of $\P^{k-1}$ (resp. of $\P^k$) is equal to $\nu(T_1,a)$ times the class of a linear subspace. So these cohomology classes do not depend on the choice of $S$. 
Kiselman showed that in general the current $S$ is not unique \cite{Kiselman}. Blel-Demailly-Mouzali gave in \cite{Blel} conditions on $T_1$ for the uniqueness of $S$. 

We consider now the situation where $T_1$ is a general positive closed $(p_1,p_1)$-current and $T_2$ is given by a submanifold $V_2$. For simplicity, we will write $T,p,V$ instead of $T_1,p_1,V_2$ and denote by $l$ the dimension of $V$. With respect to the above case, the point $a$ is replaced by the manifold $V$. We want to define a notion of tangent current to $T$ along $V$ that corresponds to the currents $S$ above. 
Let $E$ denote the normal vector bundle to $V$ in $X$ and $\overline E$ its canonical compactification. 
Denote by $A_\lambda:\overline E\to\overline E$ the map induced by the multiplication by $\lambda$ on fibers of $E$ with $\lambda\in\C^*$.  We identify $V$ with the zero section of $E$.
The tangent currents to $T$ along $V$ will be positive closed $(p,p)$-currents on $\overline E$ which are $V$-conic, i.e. invariant under the action of $A_\lambda$. The first difficulty is that when $V$ has positive dimension, in general, no neighbourhood of $V$ in $X$ is biholomorphic to a neighbourhood of $V$ in $E$. 

Let $\tau$ be a diffeomorphism between a neighbourhood of $V$ in $X$ and a neighbourhood of $V$ in $E$ whose restriction to $V$ is identity. We assume that $\tau$ is admissible in the sense that the endomorphism of $E$ induced by the differential of $\tau$ is the identity. It is not difficult to show that such maps exist, see Lemma \ref{lemma_admissible}. Here is a main result in  this paper. It is a consequence of 
Proposition \ref{prop_existence_tangent} and Theorem \ref{th_tangent_conic} below.

\begin{theorem} \label{th_main_bis}
Let $X,V,T,E,\overline E,A_\lambda$ and $\tau$ be as above. Then the family of currents $T_\lambda:=(A_\lambda)_*\tau_*(T)$ is relatively compact and any limit current, for $\lambda\to\infty$, is a positive closed $(p,p)$-current on $E$ whose trivial extension is a positive closed $(p,p)$-current on $\overline E$. Moreover, if $S$ is such a current, it is $V$-conic, i.e. invariant under $(A_\lambda)_*$, and its de Rham cohomology class in $H^{2p}(\overline E,\C)$ does not depend on the choice of $\tau$ and $S$. 
\end{theorem}  

We will see later that the result still holds and we obtain the same family of limit currents using local admissible diffeomorphisms. This softness is very useful in the analytic calculus with tangent cones and densities while the use of global admissible diffeomorphisms is convenient for cohomology calculus.

We say that $S$ is a {\it tangent current} to $T$ along $V$. Its class in the de Rham cohomology group is called the {\it total tangent class} of $T$ along $V$. Note that this notion generalizes a notion of tangent cone in the algebraic setting where $T$ is also given by a manifold, see Fulton \cite{Fulton} for details.

The cohomology ring of $\overline E$ is generated by the cohomology ring of $V$ and the tautological $(1,1)$-class on $\overline E$. Therefore, we can decompose the cohomology class of $S$ and associate to it cohomology classes of different degrees on $V$. These classes represent different parts of the tangent class of $T$ along $V$. 

Consider now arbitrary positive closed currents $T_1,T_2$ on $X$ and the tensor product $T_1\otimes T_2$ on $X\times X$. When $T_1,T_2$ are currents of integration on manifolds $V_1$ and $V_2$, the tensor product $T_1\otimes T_2$ is the current of integration on $V_1\times V_2$. Let $\Delta$ denote the diagonal of $X\times X$. We can consider the tangent currents and the total tangent class of $T_1\otimes T_2$ along $\Delta$. The normal vector bundle to $\Delta$ is canonically isomorphic to the tangent bundle of $X$ if we identify $\Delta$ with $X$. The tangent currents and the total tangent class in this case induce the {\it density currents} and the {\it total density class}  associated with $T_1$ and $T_2$.  

Assume that $p_1+p_2\leq k$ and that there is only one tangent current $S$ to $T_1\otimes T_2$ along $\Delta$. Assume also that for  $j>k-p_1-p_2$, the current
$S$ vanishes on the pull-back of $(j,j)$-forms by the canonical projection onto $\Delta$ which is canonically identified with $X$. Then we show that $S$ is the pull-back of a unique positive closed current $S^h$ of bidegree $(p_1+p_2,p_1+p_2)$ on $X$. In this case, we call $S^h$ the {\it wedge-product} of $T_1$ and $T_2$ and denote it by $T_1\curlywedge T_2$. The notion can be extended to a finite number of currents. So the density of currents extends the theory of intersection. We believe that this notion of intersection of currents, combined with the use of our theory of super-potentials, solves Demailly's problem in \cite[p.16]{Demailly2} and needs to be developed.

In Sections \ref{section_preliminar} and \ref{section_fibration}, we will recall some basic notions on positive closed currents and we give several properties that will be used in our study. 
In particular, we introduce the $\star$-norm for currents. This norm is useful for the mass estimates  of the currents $T_\lambda$ in Theorem \ref{th_main_bis} and of wedge-products of currents. We also prove the extension results needed when dealing with blow-up and with the map $\tau$ which is not holomorphic. We then introduce the notion of horizontal dimension of a current on a projective fibration and the notion of $V$-conic currents on a vector bundle over $V$. 

Tangent currents and tangent classes will be introduced in Section \ref{section_tangent}. We will prove there, a property of semi-continuity of the tangent class which is similar to the semi-continuity of Lelong number with respect to the current.  We will also give several properties which allow to compute tangent classes. In particular, we can compute such classes using  strict transforms of current and blow-up of manifolds as in Siu's results on Lelong numbers. 
Finally, the density of currents and the first properties of a new intersection theory are presented in Section \ref{section_density}. We will compare our definition with a classical notion of intersection of $(1,1)$-currents. Some applications to dynamics are given in 
\cite{DinhNguyenTruong,DinhSibony12}.

\bigskip

\noindent
{\bf Notations.} Through the paper, we denote by $\D^k$ the unit polydisc in $\C^k$ and $\lambda\D^k$ the polydisc of radius $|\lambda|$ centered at the origin of $\C^k$ for $\lambda\in\C^*$. 

If $X$ is an oriented manifold, denote by $H^*(X,\C)$ the de Rham cohomology group of $X$ and $H^*_c(X,\C)$ the de Rham cohomology group defined by forms or currents with compact support in $X$. If $V$ is a submanifold of $X$, denote by $H^*_V(X,\C)$ the de Rham cohomology group defined in the same way using only forms or currents on $X$ whose supports intersect $V$ in a compact set.

If $T$ is a closed current on $X$ denote by 
$\{T\}$ its class in $H^*(X,\C)$. When $T$ is supposed to have compact support then $\{T\}$ denotes the class of $T$ in $H^*_c(X,\C)$. If we only assume that $\supp(T)\cap V$ is compact, then $\{T\}$ denotes the class of $T$ in $H^*_V(X,\C)$. 
The current of integration on an oriented  submanifold $Y$ is denoted by $[Y]$. Its class is denoted by $\{Y\}$.

The restriction to a submanifold $V$ of smooth forms on $X$ defines a canonical morphism from $H^*_V(X,\C)$ to $H^*_c(V,\C)$; the restriction to $V$ of a class is denoted by $\{\cdot\}_{|V}$. Currents on $V$ can be canonically sent by the embedding map to currents on $X$. This induces a natural morphism from $H^*_c(V,\C)$ to $H^*_V(X,\C)$. The composition of the above two morphisms is equal to the endomorphism of the space $\oplus H^*_V(X,\C)$ induced by the cup-product with $\{V\}$. 

The group $H^{2k}_c(X,\C)$ of maximal degree is often identified with $\C$ 
via the integration of forms of maximal degree on $X$.  
If $X$ is a compact K\"ahler manifold, the groups $H^*(X,\C)$, $H^*_c(X,\C)$ and $H^*_V(X,\C)$ are equal and we identify $H^p(X,\C)$ with the direct sum of the Hodge cohomology groups $H^{q,p-q}(X,\C)$ via the Hodge decomposition. 

\section{Positive currents and spaces of test forms} \label{section_preliminar}

In this section, we recall some basic notions  on positive currents on a complex manifold and refer the reader to  Demailly \cite{Demailly3}, de Rham \cite{deRham}, Federer \cite{Federer}, H\"ormander \cite{Hormander1, Hormander2}, Siu \cite{Siu} and to \cite{DinhSibony5} for details. We will also introduce and study some spaces of test forms which are the core of the technical part of our work. They will permit, in particular, to bound the mass of the currents $T_\lambda$ in Theorem \ref{th_main_bis} and to show that their $(q,2p-q)$-components, with $q\not=p$,  converge to 0 when $\lambda$ tends to infinity. 

Let $X$ be a complex manifold of dimension $k$. A $(p,p)$-form $\theta$ on $X$ is {\it positive} if for any point in $X$ we can write $\theta$ as a finite combination of forms of type
$$(i\gamma_1\wedge\overline\gamma_1)\wedge \ldots \wedge (i\gamma_p\wedge\overline\gamma_p),$$
where  $\gamma_1,\ldots,\gamma_p$ are $(1,0)$-forms. 
A $(p,p)$-current $T$ on $X$ is {\it weakly positive} if $T\wedge \theta$ is a positive measure for any smooth positive $(k-p,k-p)$-form $\theta$. Such a current is of order 0 and real, i.e. $T=\overline T$.

A $(p,p)$-current $T$ is {\it positive} if $T\wedge\theta$ is a positive measure for any smooth weakly positive $(k-p,k-p)$-form $\theta$. Positive currents and positive forms are weakly positive.  Positivity and weak positivity are local properties. They coincide only for bidegree $(p,p)$ with $p=0,1,k-1$ or $k$.
On a chart of $X$, in the definition of (weakly) positive current,  it suffices to use only forms $\theta$ with constant coefficients.
Positive Hermitian $(1,1)$-forms on $X$ are examples of positive forms. A $(p,p)$-current $T$ is {\it strictly positive} if for a fixed smooth Hermitian $(1,1)$-form $\beta$ on $X$ we have locally $T\geq \epsilon \beta^p$, i.e. $T-\epsilon\beta^p$ is positive, for some constant $\epsilon>0$. The definition does not depend on the choice of $\beta$. 

From now on, consider a K\"ahler manifold $X$ of dimension $k$, not necessarily compact. Let  $\omega$ be a fixed K\"ahler form on $X$. It induces a K\"ahler metric on $X$ and also metrics on the vector bundles of differential forms. This  permits to define the mass-norm for currents of order 0 on $X$. If $T$ is a current of order 0 and $K$ a Borel subset of $X$, the mass of $T$ on $K$ is denoted by $\|T\|_K$ and the mass of $T$ on $X$ is denoted by $\|T\|$. If $T$ is a (weakly) positive or negative (i.e. $-T$ is positive or weakly positive) $(p,p)$-current, the above mass-norm is equivalent to the mass of  the trace measure $T\wedge \omega^{k-p}$. Then, we identify $\|T\|_K$ with the mass of $T\wedge \omega^{k-p}$ on $K$.

We introduce now some spaces of test forms and establish properties that we will use later to estimate the mass of currents. 
Fix open subsets $W_1$ and $W_2$ of $X$ with smooth boundaries such that $W_1\cap W_2$ is relatively compact in $X$. The notions below depend on the choice of $W_1$ and $W_2$. We will use later an open neighbourhood $W_2$ of a submanifold $V$ of $X$ and positive closed currents with support in $W_1$. 

\begin{definition}\rm \label{def_star}
Let $R$ be a $(1,1)$-current of order 0 on $X$ with no mass outside  $W_1\cap \overline W_2$. We define the {\it $\star$-norm} $\|R\|_\star$ of $R$ as the infimum of the constants $c\geq 0$ such that the real and imaginary parts of $R$ satisfy
$$-c(\omega +\ddc\phi)\leq \Re(R), \Im(R)\leq c(\omega +\ddc\phi)$$
for some quasi-psh function $\phi$ on $W_1$ satisfying $\ddc\phi\geq -\omega$ on $W_1$ and 
which vanishes outside $W_2$.
By convention, if such constant does not exist, the $\star$-norm of $R$ is infinite. The definition does not change if we only assume that $\phi$ is constant  outside $W_2$. 
\end{definition}

The $\star$-norm is convenient because bounding by closed positive currents permits to compute cohomologically, avoiding uniform estimates. Several concrete examples will be constructed later using blow-ups of $X$ along its submanifolds. 

Note that when $\|R\|_\star$ is finite, $R$ is absolutely continuous with respect to the positive closed $(1,1)$-current $R':=\omega+\ddc\phi$ on $W_1$. In particular, the trace measure $R\wedge\omega^{k-1}$ of $R$ is equal to the product of a bounded function with the trace measure $R'\wedge\omega^{k-1}$ of $R'$.
It is not difficult to check that $\|\cdot\|_\star$ defines a norm on the space of $(1,1)$-currents $R$ with $\|R\|_\star$ finite. This space contains the $\Cc^2$ forms with support in $W_1\cap W_2$. We will use it as a space of test forms in order to study currents with support in $W_1$. 

\begin{definition}\rm \label{def_star_bis}
Let $\Gamma$ be a form of bidegree $(1,0)$ or $(0,1)$ vanishing outside $W_1\cap \overline W_2$ with $L^2_\loc$ 
coefficients in $W_1\cap \overline W_2$. We define the $\star$-norm of $\Gamma$ by 
$\|\Gamma\|_\star:=\|i\Gamma\wedge\overline \Gamma\|_\star^{1/2}$. If $\Gamma$ is an $L^2$ 1-form  vanishing outside $W_1\cap \overline W_2$, we define $\|\Gamma\|_\star$ as the supremum of the $\star$-norms of its bidegree $(1,0)$ and bidegree $(0,1)$ components. By convention, if $\Gamma$ is a 1-current which is not given by an $L^2_\loc$ form on $W_1\cap \overline W_2$, its $\star$-norm is infinite. 
\end{definition}

\begin{remark}\rm
A version of the $\star$-norm was introduced and used by the authors for positive closed currents on compact K\"ahler manifolds \cite{DinhSibony4,DinhSibony10,DinhSibony11}, see also Vigny \cite{Vigny}. We can easily extend it to currents of bidegree $(p,p)$, $(p,0)$ or $(0,p)$. For currents $R$ of bidegree $(p,q)$ we can consider the square root of the $\star$-norm of $\Gamma\otimes \overline \Gamma$ in $X\times X$. This quantity was implicitly used in some dynamical problems, see \cite{Dinh2}.
\end{remark}

\begin{lemma} \label{lemma_def_norm}
The map $\Gamma\mapsto \|\Gamma\|_\star$ defines a semi-norm on the space of $1$-forms $\Gamma$ vanishing outside $W_1\cap \overline W_2$
such that $\|\Gamma\|_\star$ is finite. If $\Gamma_1,\Gamma_2$ are such forms then 
$$\|\Gamma_1\wedge\overline \Gamma_2\|_\star\leq \|\Gamma_1\|_\star\|\Gamma_2\|_\star.$$ 
\end{lemma}
\proof
For the first assertion, it is enough to prove the triangle inequality for the bidegree $(1,0)$-case. 
Let $\Gamma_1$ and $\Gamma_2$ be of bidegree $(1,0)$. Since the form $i(c\Gamma_1-c^{-1}\Gamma_2)\wedge (c\overline \Gamma_1-c^{-1}\overline \Gamma_2)$ is positive for any $c>0$, we have
$$i(\Gamma_1+\Gamma_2)\wedge (\overline \Gamma_1+\overline \Gamma_2)\leq (1+c^2)(i\Gamma_1\wedge \overline \Gamma_1)+(1+c^{-2})(i\Gamma_2\wedge \overline \Gamma_2).$$
Taking $c=\|\Gamma_1\|_\star^{-1/2}\|\Gamma_2\|_\star^{1/2}$,
we obtain that $\|\Gamma_1+\Gamma_2\|_\star \leq \|\Gamma_1\|_\star+\|\Gamma_2\|_\star$. 

For the second assertion, the positivity of the above form with $c=1$ implies 
$$2\Im(\Gamma_1\wedge \overline \Gamma_2)=-i\Gamma_1\wedge \overline \Gamma_2 +i\Gamma_2\wedge \overline \Gamma_1 \leq (i\Gamma_1\wedge \overline \Gamma_1)+(i\Gamma_2\wedge \overline \Gamma_2).$$
This together with similar inequalities for $\{\Gamma_1,-\Gamma_2\}$ or $\{\Gamma_1,\pm i\Gamma_2\}$ instead of $\{\Gamma_1,\Gamma_2\}$ imply 
 that $\|\Im(\Gamma_1\wedge\overline \Gamma_2)\|_\star\leq \|\Gamma_1\|_\star\|\Gamma_2\|_\star$ and a similar inequality for
$\Im(\Gamma_1\wedge\overline \Gamma_2)$.
This completes the proof of the lemma. 
\endproof 

Note that when $X$ is compact and $W_1=X$, $\|\cdot\|_\star$ is in fact a norm because the mass of $\omega+\ddc \phi$ does not depend on $\phi$. Therefore, if $\|\Gamma\|_\star=0$, the mass of $\Gamma\wedge\overline\Gamma$ vanishes and hence $\Gamma=0$.  

\begin{definition} \rm \label{def_continuous_current}
A current $R$ of order 0 on $X$ is said to be {\it a quasi-continuous current or quasi-continuous form} if it vanishes outside some open set $W_R\Subset X$ and is given on that open set by a continuous form $\Theta$ which is also an $L^1$ form. The open set of points $x\in W_R$ such that $\Theta(x)\not=0$ is called {\it an essential support} of $R$.   
\end{definition}

Note that the essential support of $R$ depends on the choice of $\Theta$ and $W_R$; it is unique up to a set of zero Lebesgue measure. If $T$ is a current of order 0 and if $R$ is as above, we can define the current $T\wedge R$ on $W_R$. If this wedge-product has finite mass, we can extend it by 0 to a current on $X$ that we still denote by $T\wedge R$. The wedge-products we will consider below satisfy this property. We will estimate their mass not in terms of the $L^\infty$ norm of $R$, as is usual.

\begin{lemma} \label{lemma_mass_star}
Let $T$ be a positive closed $(p,p)$-current on $X$ such that $\supp(T)\subset W_1$. Let $R$ be a quasi-continuous form of bidegree $(1,1)$ vanishing outside $W_1\cap W_2$ and with finite $\star$-norm.  Then 
there is a constant $c>0$ independent of $T,R$ and 
a positive closed $(p+1,p+1)$-current $T'$ on $X$ such that for any open set 
$W$ which contains $\overline W_2$
$$-T'\leq \Re(T\wedge R),\Im(T\wedge R)\leq T'\quad \mbox{and} \quad \|T'\|_W=c\|T\|_W\|R\|_\star.$$ 
\end{lemma}
\proof
We can assume that $R$ is a real current such that $\|R\|_\star=1/2$ and that there is a quasi-psh function $\phi$ on $W_1$ which vanishes outside $W_2$ and satisfies $\ddc\phi\geq -\omega$ and $-\ddc\phi-\omega\leq R\leq\ddc\phi+\omega$ on $W_1$. Since $R$ vanishes outside $W_2$, we can assume that $W_R\subset W_1\cap W_2$. 
Let $\chi_n$ be a sequence of smooth functions, with compact support in $W_R$, with $0\leq \chi_n\leq 1$ and which increases to the characteristic function of $W_R$. Define $R_n:=\chi_n R$. We still have $-\ddc\phi-\omega\leq R_n\leq\ddc\phi+\omega$.

In order to regularize $\phi$, we apply Demailly's method which uses local convolution operators, see \cite{Demailly3}. These operators act also on smooth forms $R_n$ and do not change $R_n$ too much. We can find a smooth function $\phi_n$ on an open subset $W_{1,n}$ of $W_1$ which vanishes out of an open set $W_{2,n}\supset W_2$ and  such that $\ddc\phi_n\geq -\omega$  and $-c(\ddc\phi_n+\omega)\leq R_n\leq c(\ddc\phi_n+\omega)$, where $c>0$ is a constant independent of $\phi,R$ and $n$. The constant $c$ takes into account the loss of positivity in the regularization procedure.  We can choose $W_{1,n}$ increasing to $W_1$ with $W_{1,n}\supset \supp(T)$ and $W_{2,n}$ decreasing to $W_2$.

Define $T_n:=cT\wedge (\ddc\phi_n+\omega)$. For $n$ large enough, this current is supported by $\supp(T)\subset W_{1,n}$ and the restriction of $\phi_n$ to $\supp(T)$ vanishes outside a compact subset of $W$. This and Stokes' formula imply that 
$$\|T_n\|_W =\langle T_n,\omega^{k-p-1}\rangle_W=c\langle T,\omega^{k-p}\rangle_W=c\|T\|_W.$$
In particular, the mass of $T_n$ is locally bounded uniformly on $n$. Extracting a subsequence, we can assume that $T_n$ converges to a current $T'$. Clearly, $T'$ satisfies the lemma. 
\endproof

\begin{definition}\rm \label{def_negligible}
Let $(R_\lambda)$ be a family of $q$-currents on $X$ with $\lambda\in\C$ and $|\lambda|\geq 1$. Assume that they have no mass outside $W_1\cap\overline W_2$. We say that $(R_\lambda)$ is  {\it $\star$-negligible} if it can be written as a finite sum of families of $q$-currents of type 
$$\Gamma_\lambda^1\wedge \ldots\wedge \Gamma_\lambda^q$$
where for each index $j$, the $\Gamma_\lambda^j$ are quasi-continuous forms of the same bidegree $(1,0)$ or $(0,1)$ with $\star$-norms bounded uniformly on $\lambda$ and such that one of the following properties holds
\begin{enumerate}
\item[(a)] The number of $(1,0)$-forms is not equal to the number of $(0,1)$-forms;
\item[(b)] For some index $j$, the $\star$-norm of $\Gamma_\lambda^j$ tends to 0 as $\lambda$ tends to infinity;
\item[(c)] For at least $q-1$ indices $j$, we can write $\Gamma^j_\lambda=h_\lambda^jS^j$ where $S^j$ is a quasi-continuous form independent of $\lambda$ with finite $\star$-norm and $h_\lambda^j$ are quasi-continuous functions, bounded uniformly on $\lambda$, whose essential support converges to the empty set as  $\lambda$ tends to infinity. 
\end{enumerate}
\end{definition}

Here, we say that a family of open sets $(U_\lambda)$ converges to the empty set if the characteristic function $\ind_{U_\lambda}$ converges pointwise to 0.  
The following lemma justifies the introduction of $\star$-negligible families of currents.

\begin{lemma} \label{lemma_negligible}
Let $T$ be a positive closed $(p,p)$-current on $X$ with support in $W_1$. Let $(R_\lambda)$ be a $\star$-negligible family of $(2k-2p)$-forms. Then the mass of $T\wedge R_\lambda$ converges to $0$ when $\lambda$ tends to infinity.
\end{lemma}
\proof
We only have to consider the case where $R_\lambda$ is equal to the wedge-product  $\Gamma_\lambda^1\wedge \ldots\wedge \Gamma_\lambda^{2k-2p}$ as in Definition \ref{def_negligible}. If it satisfies property (a) in that definition, then for bidegree reason, we have $T\wedge R_\lambda=0$. So we can, without loss of generality, assume that $\Gamma^j_\lambda$ is of bidegree $(1,0)$ when $j\leq k-p$ and of bidegree $(0,1)$ otherwise. 
Denote for simplicity $\Lambda_\lambda^j:=\overline \Gamma^{k-p+j}_\lambda$. 

Consider now the case  (b).
We can assume that $\|\Gamma_\lambda^1\|_\star$ converges to 0. 
Observe that $\Gamma^j_\lambda\wedge \overline \Lambda^j_\lambda$ can be written in a canonical way as a linear combination with constant coefficients of the currents 
$$i\Gamma^j_\lambda\wedge \overline \Gamma^j_\lambda, \quad i\Lambda_\lambda^j\wedge \overline \Lambda_\lambda^j,\quad
i(\Gamma_\lambda^j+  \Lambda_\lambda^j)\wedge (\overline \Gamma_\lambda^j+\overline \Lambda_\lambda^j)\quad \mbox{and} \quad  i(\Gamma_\lambda^j+i \Lambda_\lambda^j)\wedge (\overline \Gamma_\lambda^j+\overline {i\Lambda}_\lambda^j).$$ 
Therefore, we can assume that $\Lambda_\lambda^j=\Gamma_\lambda^j$ for $j\geq 2$. Define 
$$T_\lambda:=T\wedge (i\Gamma_\lambda^2\wedge \overline \Gamma_\lambda^2)\wedge\ldots\wedge  (i\Gamma_\lambda^{k-p}\wedge \overline \Gamma_\lambda^{k-p}).$$
This is a positive current. 

If we apply Lemma \ref{lemma_mass_star} inductively $k-p-1$ times to $R:=\Gamma_\lambda^j\wedge \overline \Gamma_\lambda^j$ we get a positive closed $(k-1,k-1)$-current $T'_\lambda$ of bounded mass such that $T_\lambda\leq T'_\lambda$.
Finally, Cauchy-Schwarz's inequality implies that
$$\|T_\lambda\wedge \Gamma^1_\lambda\wedge \overline \Lambda^1_\lambda\|\leq \|T_\lambda\wedge (i\Gamma_\lambda^1\wedge \overline \Gamma_\lambda^1)\|^{1/2} \|T_\lambda \wedge (i\Lambda_\lambda^1\wedge \overline \Lambda_\lambda^1)\|^{1/2}.$$
Applying again Lemma \ref{lemma_mass_star} to $T'_\lambda$ instead of $T$ and to $R=i\Gamma^1_\lambda\wedge\overline \Gamma^1_\lambda$ or $R=i\Lambda_\lambda^1\wedge \overline \Lambda_\lambda^1$ gives the result.

Assume now that condition (c) is satisfied. We can assume that it holds for $j\not=k-p+1$ and that the functions $h^j_\lambda$ in this condition are the characteristic functions of open sets $W_\lambda$ which converge to the empty set. As in the last case, we reduce the problem to the case where $\Gamma^j_\lambda=\Lambda^j_\lambda$ for $2\leq j\leq k-p$; these currents are also equal to $S^j$ restricted to $W_\lambda$. Consider the positive current
$$\widetilde T:=T\wedge (iS^2\wedge \overline S^2)\wedge\ldots\wedge  (iS^{k-p}\wedge \overline S^{k-p}).$$
We obtain as above
$$\|\widetilde T\wedge \Gamma^1_\lambda\wedge \overline \Lambda^1_\lambda\|_{W_\lambda}\leq \|\widetilde T\wedge (iS^1\wedge \overline S^1)\|_{W_\lambda}^{1/2} \|\widetilde T \wedge (i\Lambda_\lambda^1\wedge \overline \Lambda_\lambda^1)\|^{1/2}.$$
The first factor in the right hand side tends to 0 since this is the mass of a fixed current on open sets which converge to the empty set. The second factor is bounded according to 
Lemma \ref{lemma_mass_star}. The result follows.
\endproof

\begin{lemma} \label{lemma_star}
Let $M$ be a positive constant.
Let $\varphi_1,\varphi_2$ be quasi-psh functions on $W_1$ which are constant outside $W_2$ and satisfy $\ddc\varphi_1\geq -M\omega$ and $\ddc \varphi_2\geq -M\omega$. Then $\phi:=\log(e^{\varphi_1}+e^{\varphi_2})$ is a quasi-psh function on $W_1$. It is constant outside $W_2$ and satisfies $\ddc\phi\geq -M\omega$. Moreover, the $\star$-norm of
$$\Gamma:={e^{{1\over 2}(\varphi_1+\varphi_2)}\over 
e^{\varphi_1}+e^{\varphi_2}}(\partial\varphi_1-\partial\varphi_2)$$
is bounded by $\sqrt{8\pi M}$.
\end{lemma}
\proof
Clearly, $\phi$ is constant outside $W_2$. Define $\chi(t):=\log(1+e^t)$. We have $0\leq\chi'(t)\leq 1$ and $\chi''(t)\geq 0$. Therefore, if $t:=\varphi_1-\varphi_2$ we have
\begin{eqnarray*}
\ddc\phi & = & \ddc(\chi(\varphi_1-\varphi_2))+\ddc\varphi_2\\
& = & \chi'(t)(\ddc\varphi_1-\ddc\varphi_2)+{1\over 2\pi}\chi''(t)(i\partial t\wedge\overline\partial t)+\ddc\varphi_2 \\
& \geq & \chi'(t)\ddc\varphi_1+(1-\chi'(t))\ddc\varphi_2.
\end{eqnarray*}
Recall that $\ddc={i\over 2\pi}\ddbar$. 
Hence, $\phi$ is quasi-psh on $W_1$  and $\ddc\phi\geq -M\omega$. We deduce that the $\star$-norm of $\ddc\phi$ is bounded by $2M$.

A direct computation as above gives that $i\ddbar \phi-i\Gamma\wedge\overline \Gamma$ is equal to
$${e^{\varphi_1}\over e^{\varphi_1}+e^{\varphi_2}}(i\ddbar\varphi_1)+{e^{\varphi_2} \over 
 e^{\varphi_1}+e^{\varphi_2}}(i\ddbar\varphi_2).$$
The $\star$-norm of the last sum is bounded by  $4\pi M$ because $i\ddbar\varphi_1$ and $i\ddbar\varphi_2$ satisfy the same property. This gives the last estimate in the lemma.
\endproof

We will use test forms with finite $\star$-norms. The results obtained above permit to bound some integrals without knowing the $L^\infty$-norm of test forms which is not controlled by the $\star$-norm. Therefore, specific test forms with bounded $\star$-norms can be used to study singularities of currents.
We now describe a situation that we will consider in the next section, in particular, when dealing with a blow-up.

Let $V$ and $V'$ be  submanifolds of $X$ of dimension $l$ and $l'$ respectively such that $V'\subset V\subset W_2$. So $V\cap W_1$ is relatively compact in $X$. 
In $W_1\cap W_2$ consider a chart  which is identified with the polydisc $2\D^k$ on which $V$ and $V'$ are equal respectively to $2\D^l\times\{0\}$ and $2\D^{l'}\times \{0\}$. We will use in this polydisc the standard  coordinates $x=(x^1,x^2,x^3)$ with $x^1=(x_1,\ldots,x_{l'})$, $x^2=(x_{l'+1},\ldots,x_l)$ and $x^3=(x_{l+1},\ldots,x_k)$. 
In particular, $2\D^l\times\{0\}$ is the set of points $(\cdot,\cdot,0)$ in these coordinates. 
The following lemmas introduce useful families of test forms.

\begin{lemma} \label{lemma_star_loc_2}
Let $\Gamma_{m,j}$ denote the $(1,0)$-form supported by $\D^k$ and defined by
$$\Gamma_{m,j}:={x_jdx_m-x_mdx_j\over \|x^3\|^2} \quad \mbox{for}\quad l+1\leq m,j\leq k.$$
Then $\Gamma_{m,j}$ has finite $\star$-norm.
\end{lemma}

We first introduce some notations. 
Let $\sigma:\widehat X\to X$ denote the blow-up of $X$ along $V$ and define $\widehat V:=\sigma^{-1}(V)$. 
By Blanchard's theorem \cite{Blanchard}, if $U$ is an open subset of $\widehat X$ such that $U\cap \widehat V$ is relatively compact in $\widehat X$, then $U$ is a 
K\"ahler manifold. Let $\widehat\omega$ be a K\"ahler form on a neighbourhood of $\widehat W_1:=\sigma^{-1}(W_1)$. We can choose $\widehat\omega$ so that $\sigma_*(\widehat\omega)$ is equal to a constant times $\omega$ outside $W_2$. 
The current $\sigma_*(\widehat\omega)$ is positive closed and has positive Lelong number along $V$ (if $V$ is a hypersurface then $\sigma=\id$; we replace $\sigma_*(\widehat\omega)$ by $\omega+[V]$). 
Multiplying $\widehat\omega$ with a constant allows  to assume that the Lelong number of $\sigma_*(\widehat\omega)$ along $V$ is equal to 1 or equivalently, if $\widehat V:=\sigma^{-1}(V)$ is the exceptional hypersurface then $\sigma^*(\sigma_*(\widehat\omega))=\widehat\omega +[\widehat V].$

Since $\sigma_*(\widehat\omega)$ is smooth outside $V$, we can find a negative quasi-psh function $\varphi$ on $W_1$ which vanishes outside $W_2$ and such that $\ddc\varphi-\sigma_*(\widehat\omega)$ is a smooth form.
Fix a constant $c_0>1$ large enough such that $\ddc\varphi\geq \sigma_*(\widehat\omega)-(c_0-1)\omega$ and 
define $\alpha:=\ddc\varphi+c_0\omega$. This form is larger than $\sigma_*(\widehat\omega)+\omega$ and its restriction to $W_1$ has finite $\star$-norm.

We cover $\sigma^{-1}(\D^k)$ with $k-l$ charts. We describe only one of them. The other ones are obtained 
by permuting the coordinates. The chart we consider is denoted by $\widehat D$ and is given with local coordinates $z=(z_1,\ldots,z_k)$ with $|z_i|<1$ for $i\leq l$ and $|z_i|<2$ otherwise and such that
$$\sigma(z)=(z_1,\ldots,z_l,z_{l+1}z_k,\ldots,z_{k-1}z_k,z_k).$$
On this chart, $\widehat V$ is equal to $\{z_k=0\}$. Since $\ddc(\varphi\circ\sigma)-[\widehat V]$ is a smooth form and $\ddc\log|z_k|=[\widehat V]$, the function $\varphi\circ\sigma-\log|z_k|$ is smooth. We deduce  that  $\varphi-\log\|x^3\|$ is a bounded function.

\medskip\noindent
{\bf Proof of Lemma \ref{lemma_star_loc_2}.} When $V$ is a hypersurface, i.e. $l=k-1$, we have $\Gamma_{m,j}=0$. 
Consider the higher codimension case. Observe that 
 $2i\ddbar \log\|x^3\|$ is equal to the sum of 
$i\Gamma_{m,j}\wedge \overline\Gamma_{m,j}$. So we only have to check that the form 
 $i\ddbar \log\|x^3\|$
restricted to $\D^k$ has finite $\star$-norm.

Using the local coordinates $z$ introduced above, we see that 
$\sigma^*(i\ddbar\log\|x^3\|)$ is bounded by $2\pi [\widehat V]$ plus a smooth $(1,1)$-form. It follows that $i\ddbar\log\|x^3\|$ is bounded by a constant times the form $\alpha$ introduced above.  Therefore, $i\ddbar\log\|x^3\|$ has finite $\star$-norm.
\hfill $\square$

\medskip

Denote by $A_\lambda$ the map $(x^1,x^2,x^3)\mapsto (x^1,x^2,\lambda x^3)$ for $\lambda\in \C^*$. We will be concerned with $|\lambda|\to\infty$.  Therefore, in what follows, we assume that $|\lambda|\geq 1$.

\begin{lemma} \label{lemma_star_loc_1}
Let $R$ (resp. $\Gamma$)  be a quasi-continuous form essentially supported in $\D^k$ and of bidegree $(1,1)$ (resp. $(1,0)$ or $(0,1)$). Assume that their coefficients have modulus smaller or equal to $1$. Then the forms $(A_\lambda)^*(R)$ and $(A_\lambda)^*(\Gamma)$ on $\D^l\times \lambda^{-1}\D^{k-l}$ have $\star$-norms bounded  by a constant independent of $R,\Gamma$ and $\lambda$. 
\end{lemma}
\proof
Observe that the estimate on $\|(A_\lambda)^*(\Gamma)\|_\star$ can be deduced from the estimate on $\|(A_\lambda)^*(R)\|_\star$ applied to $R:=\Gamma\wedge\overline \Gamma$. So it is enough to bound the $\star$-norm of $(A_\lambda)^*(R)$. If $R$ does not contain terms with $dx_j$ or $d\overline x_j$ with $j\geq l+1$, then $(A_\lambda)^*(R)$ has bounded coefficients and its $\star$-norm is clearly bounded. Moreover, since we can bound the real and imaginary parts of $dx_m\wedge d\overline x_j$ by $idx_m\wedge d\overline x_m+idx_j\wedge d\overline x_j$, we only have to consider the case where $R=i\ddbar\|x^3\|^2$. 

We construct now a function $\phi$ satisfying estimates as in Definition \ref{def_star}. Define $s:=\log|\lambda|$. 
Recall that $\varphi-\log\|x^3\|$ is a bounded function, where the function $\varphi$ was defined above. Fix a constant $A\geq 1$ large enough such that $-A\leq \varphi-\log\|x^3\|\leq A$. We only have to consider the case where $s$ is large enough, e.g. $s\geq 3A$. 
Observe that since $A$ is large enough, we have $\varphi\leq -s+2A$ on $\D^l\times \lambda^{-1}\D^{k-l}$.

Let $\chi$ be a convex increasing function on $\R$ such that $\chi(t)=t$ for $t\geq -s+3A$, $0\leq \chi'\leq 1$ everywhere and $\chi''(t)= e^{2t+2s-5A}$ for $t\leq -s+2A$. 
Define $\phi:=c^{-1}\chi\circ\varphi$ for a fixed constant $c\gg c_0$ large enough. 
It is clear that $\phi$ vanishes outside $W_2$. 
A direct computation gives
$$i\ddbar\phi=c^{-1}\chi''(\varphi) i\partial\varphi\wedge \dbar\varphi + c^{-1} \chi'(\varphi) i\ddbar\varphi.$$
The first term in the last sum is positive. The second one is bounded below by $-2\pi c^{-1}c_0\omega$. Therefore, $i\ddbar\phi\geq -\omega$ and $\phi$ is quasi-psh. 

We prove now that $e^{2s}i\ddbar \|x^3\|^2\leq c^3(i\ddbar\phi+\alpha)$ on the open subset of $\D^k$ where  
$\varphi <-s+2A$. This property implies that the form $(A_\lambda)^*(i\ddbar\|x^3\|^2)$ on $\D^k\times \lambda^{-1}\D^{k-l}$ has bounded $\star$-norm because it is equal to $e^{2s}i\ddbar \|x^3\|^2$.
The idea is to pull-back the forms by $\sigma$ and check the inequality in the chart $\widehat D$ that we have described. 

Define $\widehat\varphi:=\varphi\circ \sigma$ and $\widehat\phi:=\phi\circ \sigma=c^{-1}\chi\circ \widehat\varphi$. 
Since $\widehat \varphi-\log|z_k|$ is a smooth function, the form $\gamma:=\partial(\widehat\varphi-\log|z_k|)$ is smooth.
Recall also that $i\ddbar\widehat\varphi\geq -2\pi c_0\sigma^*(\omega)$ and $\widehat\varphi\geq \log|z_k|-A$. Therefore, when $\widehat\varphi(z)<-s+2A$, a direct computation as above gives
\begin{eqnarray*}
c^3(i\ddbar\widehat\phi+\sigma^*(\alpha)) &\geq& c^2\chi''(\widehat\varphi)   i\partial\widehat\varphi\wedge \dbar\widehat\varphi +
c^2\chi'(\widehat\varphi)i\ddbar\widehat\varphi+c^3\sigma^*(\alpha)\\
& \geq & c^2e^{2s-7A}|z_k|^2  i\partial\widehat\varphi\wedge \dbar\widehat\varphi -2c^2\pi c_0\sigma^*(\omega)+c^3(\widehat \omega+\sigma^*(\omega))\\
& \geq & c^2e^{2s-7A}|z_k|^2i(z_k^{-1}dz_k+\gamma)\wedge(\overline z_k^{-1} d\overline z_k+\overline\gamma)+c^3\widehat\omega.
\end{eqnarray*}
We also have
\begin{eqnarray*}
2i(z_k^{-1}dz_k+\gamma)\wedge (\overline z_k^{-1}d\overline z_k + \overline\gamma)
& = &    i(z_k^{-1}dz_k+2\gamma)\wedge (\overline z_k^{-1}d\overline z_k +2 \overline\gamma)\\
& & +|z_k|^{-2} idz_k\wedge d\overline z_k -2i\gamma\wedge\overline\gamma\\
& \geq &  |z_k|^{-2} idz_k\wedge d\overline z_k -2i\gamma\wedge\overline\gamma.
\end{eqnarray*}
Since $e^{2s-7A}|z_k|^2\gamma\wedge\overline\gamma$ is a bounded form on the considered domain and because the constant $c$ is large enough, we deduce from the inequalities above that
$$c^3(i\ddbar\widehat\phi+\sigma^*(\alpha))\geq ce^{2s}idz_k\wedge d\overline z_k+c\widehat\omega.$$ 

On the other hand, one can find bounded forms $\theta_i$ on $\D^l\times\lambda^{-1}\D^{k-l}$ such that
$$\sigma^*(e^{2s}\ddc\|x^3\|^2)  =   e^{2s}\ddc|z_k|^2+e^sdz_k\wedge\theta_1+e^sd\overline z_k\wedge\theta_2+\theta_3.$$
Cauchy-Schwarz's inequality implies that the last sum is bounded above by $2e^{2s}\ddc|z_k|^2+\theta_4$ for some bounded 
form $\theta_4$. We conclude that
$$\sigma^*(e^{2s}\ddc\|x^3\|^2)\leq c^3(i\ddbar\widehat\phi+\sigma^*(\alpha))$$
on $\widehat D\cap\sigma^{-1}(\D^l\times \lambda^{-1}\D^{k-l})$. Hence,
$$e^{2s}\ddc\|x^3\|^2\leq c^3(i\ddbar\phi+\alpha)$$
on $\D^l\times \lambda^{-1}\D^{k-l}$. This completes the proof of the lemma.
\endproof

\begin{lemma} \label{lemma_star_loc_3}
Let $\Gamma'_{m,j}$ be the $(1,0)$-form on $\C^k$ 
given by 
$$\Gamma'_{m,j}:={x_jdx_m-x_mdx_j\over \|x^2\|^2+\|x^3\|^2} \quad \mbox{for}\quad l'+1\leq m,j\leq k.$$ 
Then  the restriction of  
$(A_\lambda)^*(\Gamma'_{m,j})$ to $\D^k$ has $\star$-norm bounded by a constant which does not depend on $\lambda,m$ and $j$ with $|\lambda|\geq 1$.
\end{lemma}
\proof
We have  to bound the $\star$-norm of $(A_\lambda)^*(i\Gamma'_{m,j}\wedge\overline \Gamma'_{m,j})$. Observe that the form $i\Gamma'_{m,j}\wedge\overline \Gamma'_{m,j}$ is positive and bounded by  the form $R:=i\ddbar\log(\|x^2\|^2+\|x^3\|^2)$. So it is enough to bound the $\star$-norm of the restriction of
$(A_\lambda)^*(R)=i\ddbar \log(\|x^2\|^2+|\lambda|^2\|x^3\|^2)$ to $\D^k$. Write $\widetilde\varphi_1:=\log[(|\lambda|^2-1)\|x^3\|^2]$ and $\widetilde\varphi_2:=\log(\|x^2\|^2+\|x^3\|^2)$. We have $(A_\lambda)^*(R)=i\ddbar\log(e^{\widetilde\varphi_1}+e^{\widetilde\varphi_2})$.

With notations as in Lemma \ref{lemma_star_loc_2}, we can write $i\ddbar\widetilde\varphi_1$ as a finite combination of 
$\Gamma_{m,j}\wedge\overline\Gamma_{m,j}$. Therefore, $i\ddbar\widetilde\varphi_1$ has bounded $\star$-norm. 
The same arguments applied to $V'$ instead of $V$ imply that the $\star$-norm of  $i\ddbar\widetilde\varphi_2$ is bounded.
Recall that the function $\varphi$ was defined above using the blow-up $\sigma:\widehat X\to X$ of $X$ along $V$. 
Let $\varphi'$ be the function obtained in the same way by replacing $V$ with $V'$.
Define also $\varphi_1:=2\varphi+\log(|\lambda|^2-1)$ and $\varphi_2:=2\varphi'$. 
Observe that $i\ddbar\varphi_1$ and $i\ddbar \varphi_2$ restricted to $\D^k$ have $\star$-norms bounded independently of $\lambda$. 

Using the coordinates $z$ on $\widehat D$ introduced above, we see that $(\varphi_1-\widetilde\varphi_1)\circ\sigma$ 
is the potential of a smooth form. So it is a smooth function. We deduce that
$\varphi_1-\widetilde\varphi_1$ is a bounded function and $\partial\varphi_1-\partial\widetilde\varphi_1$ is the push-forward by $\sigma$ of a combination of $dz_j$ with bounded coefficients. Therefore, $\partial\varphi_1-\partial\widetilde\varphi_1$ is equal to the sum of a bounded form and a combination with bounded coefficients of $\Gamma_{m,j}$. We deduce that this form has bounded $\star$-norm. In the same way, we obtain that $\varphi_2-\widetilde\varphi_2$ is a bounded function and $\partial\varphi_2-\partial\widetilde\varphi_2$ has bounded $\star$-norm. These functions and forms do not depend on $\lambda$.

Define 
$$\Gamma:={e^{{1\over 2}(\varphi_1+\varphi_2)}\over 
e^{\varphi_1}+e^{\varphi_2}}(\partial\varphi_1-\partial\varphi_2) \quad \mbox{and} \quad 
\widetilde \Gamma:={e^{{1\over 2}(\widetilde\varphi_1+\widetilde\varphi_2)}\over 
e^{\widetilde\varphi_1}+e^{\widetilde\varphi_2}}(\partial\widetilde\varphi_1-\partial\widetilde\varphi_2).$$
The coefficients involving the exponential in the last line are smaller than 1 since the exponential function is convex.  
We deduce from the above  properties of $\partial\varphi_j-\partial\widetilde\varphi_j$ 
that $\widetilde \Gamma$ is equal to a $\star$-bounded form plus the form
$$\widetilde \Gamma':={e^{{1\over 2}(\widetilde\varphi_1+\widetilde\varphi_2)}\over 
e^{\widetilde\varphi_1}+e^{\widetilde\varphi_2}}(\partial\varphi_1-\partial\varphi_2).$$
Moreover, $\widetilde \Gamma'$ is equal to a bounded function times $\Gamma$ which is $\star$-bounded 
according to Lemma \ref{lemma_star}. It follows that $\|\widetilde \Gamma\|_\star$ is also bounded. 

A computation as in Lemma \ref{lemma_star} shows that
$$(A_\lambda)^*(R)=i\widetilde \Gamma\wedge\overline{\widetilde \Gamma} + {e^{\widetilde\varphi_1}\over e^{\widetilde\varphi_1}+e^{\widetilde\varphi_2}}(i\ddbar\widetilde\varphi_1)+{e^{\widetilde\varphi_2} \over 
 e^{\widetilde\varphi_1}+e^{\widetilde\varphi_2}}(i\ddbar\widetilde\varphi_2).$$ 
By Lemma \ref{lemma_star_loc_2},  $i\ddbar\widetilde\varphi_1$ and  $i\ddbar\widetilde\varphi_2$ have bounded $\star$-norms. It follows that  $(A_\lambda)^*(R)$ has also a bounded $\star$-norm. This completes the proof of the lemma.
\endproof

\begin{lemma} \label{lemma_star_exact}
Let $R$ be a smooth $q$-form with compact support in $\D^k$. Let $t\in\C$ be a fixed constant such that $|t|\geq 1$. Assume that $V$ is a hypersurface, i.e. $l=k-1$ and $\sigma=\id$. Then there are smooth $(q-1)$-forms $\Theta_\lambda$ with compact support in $\D^k$ such that the family of forms 
$$(A_{t\lambda})^*(R)-(A_\lambda)^*(R)-d\Theta_\lambda$$ 
on $\D^l\times \lambda^{-1}\D^{k-l}$ is $\star$-negligible.  
\end{lemma}
\proof
The map $A_\lambda$ is given by $(x_1,\ldots,x_{k-1},x_k)\mapsto (x_1,\ldots,x_{k-1},\lambda x_k)$. 
Solving the d-equation on the complex lines where $x_1,\ldots,x_{k-1}$ are constant, we obtain a smooth $(q-1)$-form $\Theta$ with compact support in $\D^k$ such that $A_{t}^*(R)-R-d\Theta$
does not contain terms with $dx_k\wedge d\overline x_k$. Define $\Theta_\lambda:=A_\lambda^*(\Theta)$. We have
$$A_{t\lambda}^*(R)-A_\lambda^*(R)-d\Theta_\lambda=A_\lambda^*\big(A_{t}^*(R)-R-d\Theta\big).$$
So these forms do not contain terms with  $dx_k\wedge d\overline x_k$. 

Finally, observe that $dx_j$ and $d\overline x_j$ are invariant under the actions of $A_\lambda^*$ and $A_{t\lambda}^*$ for $j\leq k-1$. By Lemma \ref{lemma_star_loc_1}, they have finite $\star$-norms. Therefore, the family of forms in the lemma satisfies the property (c) in Definition \ref{def_negligible}.
\endproof

Let $x':=(x^1,x^2)$ and $x'':=x^3$.
For $m\geq 0$, denote by $O^*(\|x''\|^m)$ a function (resp. a 1-form) which is continuous outside $V$ and is equal to (resp. whose coefficients are equal to) $O(\|x''\|^m)$ when $x''\to 0$.
Recall that a function is equal to $O(\|x''\|^m)$ if its modulus is bounded by a constant times $\|x''\|^m$. We will use functions and forms depending on a parameter $\lambda$ and we always assume that the constant is independent of $\lambda$.  
Denote also by  $O^{**}(\|x''\|^m)$ the sum of a 1-form with $O^*(\|x''\|^m)$ coefficients and a combination of $dx'', d\overline x''$ with $O^*(\|x''\|^{m-1})$ coefficients for $m\geq 1$. 
A vector or a matrix whose coefficients satisfy the same property is denoted with the same notation.  

Consider now a bi-Lipschitz map $\tau$ from a neighbourhood $U$ of $\overline \D^k\cap V$ to another neighbourhood of 
$\overline \D^k\cap V$. We assume that $\tau$ is  
smooth outside $U\cap V$ and that its restriction to $U\cap V$ is identity. In the following expressions, we consider $x'$ and $x''$ as line matrices. 

\begin{definition} \rm \label{def_admissible_loc}
We say that $\tau$ is {\it admissible} if there 
is an $O^*(1)$ matrix $a(x)$ on $U$ such that 
$$\tau=\big(x'+x''a(x),x''\big)+O^*(\|x''\|^2)$$
and
$$d\tau(x)  =  \big(dx'+dx''a(x)+ O^*(\|x''\|), dx'' + O^{**}(\|x''\|^2)\big)$$
when $x''\to 0$.
\end{definition}

Equivalently, there are  $O^*(1)$ matrices $a_1(x)$ and $a_2(x)$ on $U$ such that 
$$\tau=\big(x^1+x^3a_1(x),x^2+x^3a_2(x),x^3\big)+O^*(\|x^3\|^2)$$
and
\begin{eqnarray*}
d\tau(x) & = & \big(dx^1+dx^3a_1(x)+ O^*(\|x^3\|), dx^2+dx^3a_2(x)+ O^*(\|x^3\|), \\
&& \qquad \qquad \qquad dx^3 +O^{**}(\|x^3\|^2)\big)
\end{eqnarray*}
when $x^3\to 0$.

Note that the differential of a smooth and admissible map is $\C$-linear at every point of $U\cap V$ but in general it does not depend holomorphically on the point. The map $\tau=(x'+x''a(x),x'')$ with $a$ holomorphic is admissible, but we will need later global admissible maps which are not necessarily holomorphic, even not smooth.

\begin{definition}\rm \label{def_principal}
Let $(R_\lambda)$ be a family of $(q,q)$-currents on $X$ with $\lambda\in\C$ and $|\lambda|\geq 1$. We say that this family is {\it $\star$-principal} if it can be written as a finite sum of families of $(q,q)$-currents of type
$$\Gamma_\lambda^1\wedge \ldots\wedge \Gamma_\lambda^{2q},$$
where the $\Gamma_\lambda^j$ are quasi-continuous forms with $\star$-norms bounded uniformly on $\lambda$ such that $q$ of them are of bidegree $(1,0)$ and $q$ are of bidegree $(0,1)$.
\end{definition} 

Let $\tau$ be an admissible map as above.
Denote by $\sigma':\widehat X'\to X$ the blow-up along $V'$. We will need some test forms which are the push-forward of smooth forms by $\sigma'$.  We have the following lemma.

\begin{lemma} \label{lemma_star_admissible}
Let $R$ be a smooth $2q$-form with compact support in $\sigma'^{-1}(\D^k)$.
Let $R'$ be the bidegree $(q,q)$ component of $R$ and define $R_\lambda:=(A_\lambda)^*\sigma'_*(R')$.
Then the family $R_\lambda$ is  $\star$-principal and the family 
$\tau^*(A_\lambda)^*\sigma'_*(R)-R_\lambda$ is $\star$-negligible. 
\end{lemma}
\proof
We cover $\sigma'^{-1}(\D^k)$ with a finite number of charts. Using a partition of unity, we can assume that $R$ is supported by one of these charts. We will only work in the chart $\widehat D'$ that we describe now. The result holds for the other ones because they are obtained from $\widehat D'$ just by using some permutations of indices. 

We have on $\widehat D'$ holomorphic coordinates $w=(w_1,\ldots,w_k)$ with $|w_i|<1$ for $i\leq l'$ and $|w_i|<2$ otherwise such that 
$$\sigma'(w)=(w_1,\ldots,w_{l'},w_{l'+1}w_k,\ldots,w_{k-1}w_k, w_k).$$
Denote by $x$ the image of $w$ by $\sigma'$. We have $w_j=x_j/x_k$ for $l'+1\leq j\leq k-1$ and $w_j=x_j$ otherwise.
We also have $\|x^2\|\lesssim |x_k|$ and $\|x^3\|\lesssim |x_k|$ on $\sigma'(\widehat D')$. This implies that 
$\|x^2\|\lesssim |\lambda||x_k|$ and $\|x^3\|\lesssim |x_k|\lesssim |\lambda|^{-1}$  on 
$(A_\lambda)^{-1}\sigma'(\widehat D')$ and on $\tau^{-1}(A_\lambda)^{-1}\sigma'(\widehat D')$ because $\tau$ is bi-Lipschitz. The later sets contain respectively  
the support of $(A_\lambda)^*\sigma'_*(R)$ and the support of $\tau^*(A_\lambda)^*\sigma'_*(R)$. 

In order to obtain the result, we first study the actions of $(A_\lambda)^*\sigma'_*$ and of $\tau^*(A_\lambda)^*\sigma'_*$ on smooth functions and on linear 1-forms. The form $R$ is built using these functions and 1-forms. 
Let $g$ be a smooth function with compact support in $\widehat D'$. If we define 
$$w_{x,\lambda}:=\sigma'^{-1}(A_\lambda(x))=(x^1,\lambda^{-1}x_k^{-1}x^2, x_k^{-1}x_{l+1},\ldots,x_k^{-1}x_{k-1},\lambda x_k),$$ 
then 
$$\tau^*(A_\lambda)^*\sigma'_*(g)(x)-(A_\lambda)^*\sigma'_*(g)(x)=
g(w_{\tau(x),\lambda})-g(w_{x,\lambda}).$$
Since $\tau$ is admissible, the above estimates on $\|x^2\|$ and $\|x^3\|$ imply that 
$\|w_{\tau(x),\lambda}-w_{x,\lambda}\|=O(\lambda^{-1})$. 
The smoothness of $g$ implies that the functions in the previous identity are uniformly bounded by a constant times $|\lambda|^{-1}$. 

Consider now the forms $\tau^*(A_\lambda)^*\sigma'_*(dw_j)$ and $\tau^*(A_\lambda)^*\sigma'_*(d\overline w_j)$. We will discuss the case of $dw_j$; the other case is treated similarly. 
Observe that since $\tau$ is admissible, for $\lambda$ large enough, the supports of the considered forms are contained in $\D^l\times 2\lambda^{-1}\D^{k-l}$. 
By Lemma \ref{lemma_star_loc_1} applied to $\lambda/2$ instead of $\lambda$, on the considered domains, bounded forms have bounded $\star$-norms and the $\star$-norm of an 
$O^{**}(\|x^3\|)$ 1-form  is of order $O(\lambda^{-1})$ as $\lambda$ tends to infinity.
We will use these properties and the admissibility of $\tau$ several times in the discussion below.

Define $w^1:=(w_1,\ldots,w_{l'})$. We have
$$\tau^*(A_\lambda)^*\sigma'_*(dw^1)=\tau^*(dx^1)=dx^1+O^{**}(\|x^3\|).$$
Since the components of $(A_\lambda)^*\sigma'_*(dw^1)=dx^1$ are bounded forms, they have bounded 
$\star$-norms. The $\star$-norm of the components of $O^{**}(\|x^3\|)$ on the considered domains tends to 0 as $\lambda$ tends to infinity.
So the $\star$-norm of $\tau^*(A_\lambda)^*\sigma'_*(dw^1) -(A_\lambda)^*\sigma'_*(dw^1)$ tends to 0.

For $j=k$,  we have 
$$\tau^*(A_\lambda)^*\sigma'_*(dw_k)=\lambda dx_k + \lambda O^{**}(\|x^3\|^2)=(A_\lambda)^*\sigma'_*(dw_k) + O^{**}(\|x^3\|).$$
As we already wrote above, the form 
$(A_\lambda)^*\sigma'_*(dw_k)=\lambda dx_k$
has bounded $\star$-norm and the $\star$-norm of $O^{**}(\|x^3\|)$ tends to 0. So the $\star$-norm of $\tau^*(A_\lambda)^*\sigma'_*(dw_k) -(A_\lambda)^*\sigma'_*(dw_k)$ tends to 0. 

Assume that $l+1\leq j\leq k-1$. We have  
$$(A_\lambda)^*\sigma'_*(dw_j)=(A_\lambda)^*\Big({x_kdx_j-x_jdx_k\over x_k^2}\Big)={x_k dx_j-x_jdx_k\over x_k^2}\cdot$$
By Lemma \ref{lemma_star_loc_2}, this form has bounded $\star$-norm.
As above, the admissibility of $\tau$ implies that
the $\star$-norm of
$\tau^* (A_\lambda)^*\sigma'_*(dw_j) -(A_\lambda)^*\sigma'_*(dw_j)$ 
tends to 0. 

Consider now the remaining case where $l'+1\leq j\leq l$. We have 
$$(A_\lambda)^*\sigma'_*(dw_j)=(A_\lambda)^*\Big({x_kdx_j-x_jdx_k\over x_k^2}\Big)={x_k dx_j-x_jdx_k\over \lambda x_k^2}\cdot$$
Since $\|x^2\|\lesssim |\lambda||x_k|$, by Lemma \ref{lemma_star_loc_3}, this form has bounded $\star$-norm.
Using the description of $d\tau$ and that $\|x^2\|\lesssim \lambda |x_k|$, we obtain that 
$$\tau^* (A_\lambda)^*\sigma'_*(dw_j)-{x_k dx_j-x_jdx_k\over \lambda x_k^2}$$
is equal to $O^{**}(\|x^3\|)$ plus a linear combination of the forms considered in the previous case with $O(\lambda^{-1})$ coefficients.
Therefore,
the $\star$-norm of $\tau^* (A_\lambda)^*\sigma'_*(dw_j) -(A_\lambda)^*\sigma'_*(dw_j)$ 
tends to 0.

We can now apply the above discussion to each component of $R$ written in $w$-coordinates.
It is easy to deduce that the family $R_\lambda$ is $\star$-principal and 
$\tau^* (A_\lambda)^*\sigma'_*(R)-R_\lambda$ is $\star$-negligible and is a sum of forms satisfying (a) or (b) in Definition \ref{def_negligible}. 
\endproof

The following lemma together with Lemmas \ref{lemma_def_norm} and \ref{lemma_mass_star} will allow to estimate forms of any degree using an induction on the degree. 

\begin{lemma} \label{lemma_star_admissible_bis}
Let $R$ be a continuous $1$-form with compact support in $\sigma'^{-1}(\D^k)$ such that $\|R\|_\infty\leq 1$. Then 
the $\star$-norm of $\tau^*(A_\lambda)^*\sigma'_*(R)$ is bounded by a constant independent of $\lambda$ and $R$.
\end{lemma}
\proof
Observe that the computation in the last lemma is valid in this case except for the estimate on $\tau^*(A_\lambda)^*\sigma'_*(g)$ when $g$ is only continuous and bounded by 1. However, we only need here that 
$\tau^*(A_\lambda)^*\sigma'_*(g)$ is bounded by 1.
For $\lambda$ large enough, $R_\lambda'$ is supported by $\D^l\times 2\lambda^{-1}\D^{k-l}$. Therefore, we easily deduce from the computation in the last lemma that 
the $\star$-norm of $\tau^*(A_\lambda)^*\sigma'_*(R)$ is bounded by a constant independent of $\lambda$ and $R$.
\endproof

Let $\tau$ be a bi-Lipschitz map from a neighbourhood $U$ of $\overline \D^k\cap V$ to another neighbourhood of 
$\overline \D^k\cap V$. We assume that $\tau$ is  
smooth outside $U\cap V$ and that its restriction to $U\cap V$ is identity. 

\begin{definition} \rm \label{def_almost_admissible_loc}
We say that $\tau$ is {\it almost-admissible} if 
$$\tau=\big(x'+O^*(\|x''\|),x''+O^*(\|x''\|^2)\big)$$
and
$$d\tau(x)  =  \big(dx'+O^{**}(\|x''\|), dx'' + O^{**}(\|x''\|^2)\big)$$
when $x''\to 0$.
\end{definition}

\begin{remark} \rm \label{rk_almost_admissible}
Let $\tau$ be an almost-admissible map as above. When $\tau$ is smooth, its differential at a point of $V$ is not necessarily $\C$-linear. 
Let $R$ be a smooth $2q$-form with compact support in $\D^k$ and let $R'$ be its component of bidegree $(q,q)$. 
Define $R_\lambda:=(A_\lambda)^*(R')$.
As in Lemma \ref{lemma_star_admissible}, we obtain that the family $\tau^*(A_\lambda)^*(R)-R_\lambda$ is $\star$-negligible.
 If $R$ is a smooth 1-form with compact support in $\D^k$ and with coefficients bounded by 1, as in Lemma \ref{lemma_star_admissible_bis}, the $\star$-norm of $\tau^*(A_\lambda)^*(R)$ is bounded by a constant independent of $\lambda$ and $R$. We can use this property together with Lemmas \ref{lemma_def_norm} and \ref{lemma_mass_star} in order to estimate forms of higher degree. 
\end{remark}

We close this section with a technical lemma that we will use in the next sections.
Let $W$ and $\widetilde W$ be K\"ahler manifolds of dimension $k$. Let $V$ and $\widetilde V$ be complex submanifolds of dimension $l$ of $W$ and $\widetilde W$ respectively. 
Consider a bi-Lipschitz map  $\tau:W\to\widetilde W$ which is smooth outside $V$, preserves the orientation and such that $\tau(V)=\widetilde V$. Denote by $\Gamma$ the graph of $\tau$ in $W\times\widetilde W$, $\Pi:W\times\widetilde W\to W$ and $\widetilde \Pi:W\times\widetilde W\to \widetilde W$ the canonical projections. 

The integration on $\Gamma$ defines a closed current $[\Gamma]$ of order 0. This can be seen using de Rham regularization theorem for currents. 
If $\theta$ is a smooth $q$-form on $W$ then $\tau_*(\theta)$ is a bounded form which is equal in the sense of currents to $\widetilde\Pi_*([\Gamma]\wedge \Pi^*(\theta))$. In particular, if $\theta$ is closed or exact, so is $\tau_*(\theta)$. It follows that $\tau$ defines a morphism $\tau_*$ from the cohomology ring $\oplus H^*_V(W,\C)$ (resp. 
$\oplus H^*_c(W,\C)$ and $\oplus H^*(W,\C)$) to the cohomology ring $\oplus H^*_{\widetilde V}(\widetilde W,\C)$ (resp. 
$\oplus H^*_c(\widetilde W,\C)$ and $\oplus H^*(\widetilde W,\C)$). The same property holds for $\tau^{-1}$ and gives a morphism $\tau^*$. 

Let $Z$ be a smooth oriented manifold of dimension $q$ and let $\rho:Z\to W$ be a Lipschitz proper map. We can define a current $\Delta_\rho$ of order 0 on $W$ 
by $\langle \Delta_\rho,\theta\rangle:=\langle Z,\rho^*(\theta)\rangle$ for smooth $q$-forms $\theta$ with compact support in $W$. We can see using the graph of $\rho$ that this current is closed. Using a regularization we obtain that $\tau_*\{\Delta_\rho\}=\{\Delta_{\tau\circ\rho}\}$. Since $\tau$ is bi-Lipschitz, we 
also obtain that $\{\Delta_\rho\}=\tau^*\{\Delta_{\tau\circ\rho}\}$. Hence, $\tau_*\circ\tau^*=\id$. Note that we have $\tau_*\{V\}=\{\widetilde V\}$ and $\tau^*\{\widetilde V\}=\{V\}$. As above, the map $\tau_{|V}:V\to \widetilde V$ 
 induces also isomorphisms $(\tau_{|V})_*$ and $(\tau_{|V})^*$ between the cohomology rings on $V$ and $\widetilde V$.

\begin{lemma} \label{lemma_bi_lipschitz}
Let $T$ be a positive closed $(p,p)$-current on $W$ without mass on $V$ such that $\supp(T)\cap V$ is compact. 
Then the current $(\tau_{|W\setminus V})_*(T)$ has finite mass on any compact subset of $\widetilde W$. 
Let $\tau_*(T)$ denote the extension of  $(\tau_{|W\setminus V})_*(T)$ by $0$ to a current 
 on $\widetilde W$. Then $\tau_*(T)$ is a closed current such that $\{\tau_*(T)\}=\tau_*\{T\}$ and 
 $\{\tau_*(T)\}_{|\widetilde V}=(\tau_{|V})_*(\{T\}_{|V})$.
\end{lemma}
\proof
Let $\theta$ be a smooth $(2k-2p)$-form with compact support on $\widetilde W$. We have
$$\langle (\tau_{|W\setminus V})_*(T),\theta\rangle = \langle T, \tau^*(\theta)\rangle_{W\setminus V}.$$
If the coefficients of $\theta$ are bounded, $\tau^*(\theta)$ satisfies the same property because $\tau$ is Lipschitz. It follows that the above integrals are bounded and then $(\tau_{|W\setminus V})_*(T)$ has  bounded mass on compact subsets of $\widetilde W$. We can extend it by 0 to a current $\tau_*(T)$ on $\widetilde W$.

We show that this current is closed. The problem concerns only the points near $V$.
Assume that $\theta=d\gamma$ with $\gamma$ smooth supported in a compact set $K$ in $\widetilde W$. We have to prove that the above integrals vanish.  Using a partition of unity, we can assume that $\tau^{-1}(K)$ is contained in the chart $\D^k$ as above. Let $\chi$ be a function on $\D^l\times \C^{k-l}$ which vanishes in a neighbourhood of $V=\D^l\times \{0\}$ and is equal to 1 outside $\D^k$. Since $T$ is closed, the considered integrals are equal to
$$\lim_{\lambda\to \infty}\big\langle T,(\chi\circ A_\lambda) d(\tau^*(\gamma))\big\rangle= \lim_{\lambda\to \infty} - \big\langle T, (A_\lambda)^*(d\chi)\wedge \tau^*(\gamma)\big\rangle.$$

We show that $\langle T, (A_\lambda)^*(\partial \chi)\wedge \tau^*(\gamma)\rangle$ tends to 0
 as $\lambda\to \infty$. This and an analogous property with $\dbar\chi$ instead of $\partial\chi$ give the result. We only have to consider the bidegree $(k-p-1,k-p)$ part of $\tau^*(\gamma)$ because $T$ is of bidegree $(p,p)$. Since this is a bounded form, we can write it as a finite combination of forms of type $\beta\wedge\Theta$, where $\beta$ is a $(0,1)$-form and $\Theta$ a positive $(k-p-1,k-p-1)$-form, both are bounded and smooth outside $V$. Without loss of generality,  we can replace $\tau^*(\gamma)$ by $\beta\wedge\Theta$. 

Define $\Gamma_\lambda:=(A_\lambda)^*(\partial\chi)$. We obtain from the Cauchy-Schwarz's inequality
$$|\langle T,\Gamma_\lambda\wedge\beta \wedge \Theta\rangle|^2\leq \langle T,i\overline\beta\wedge\beta\wedge \Theta\rangle_{\D^l\times \lambda^{-1}\D^{k-l}} \langle T,i \Gamma_\lambda\wedge\overline \Gamma_\lambda\wedge \Theta\rangle.$$
The first factor in the right hand side of the last inequality tends to 0 since $T$ has no mass on $V$. The second one is bounded according to Lemmas \ref{lemma_mass_star} and \ref{lemma_star_loc_1}. This implies that $\tau_*(T)$ is closed.

We prove now the first identity in the lemma. Let $\theta$ be a smooth closed $(2k-2p)$-form such that $\supp(\theta)\cap \supp(\tau_*(T))$ is compact. We have seen that $\tau^*(\theta)$ is a closed current. It is enough to check that $\{\tau_*(T)\}\smallsmile \{\theta\} = \{T\}\smallsmile \{\tau^*(\theta)\}$. Since $\tau$ is smooth outside $V$, by definition of $\tau_*(T)$, we have
$$\{\tau_*(T)\}\smallsmile \{\theta\} =\langle \tau_*(T),\theta\rangle = \langle T,\tau^*(\theta)\rangle_{W\setminus V}.$$
By de Rham's regularization theorem, there exists a sequence of smooth closed $(2k-2p)$-forms $\theta_n$ supported in a fixed open subset $W'$ of $W$ such that $\theta_n\to \tau^*(\theta)$ and $W'\cap \supp(T)$ is relatively compact in $W$. 
These forms are obtained from $\tau^*(\theta)$ by convolution with diffeomorphisms of $W$ which approximate the identity. Thus, 
$\tau^*(\theta)$ is smooth outside $V$ and bounded on $W$, the forms $\theta_n$ are bounded uniformly on $n$ and converge locally uniformly to $\tau^*(\theta)$ on $W\setminus V$. Finally, since $T$ has no mass on $V$, we have
$$\{T\}\smallsmile \{\tau^*(\theta)\}= \lim_{n\to\infty}\{T\}\smallsmile \{\theta_n\} = \lim_{n\to\infty} \langle T,\theta_n\rangle = \langle T,\tau^*(\theta)\rangle_{W\setminus V}.$$
Hence, $\{\tau_*(T)\}=\tau_*\{T\}$.

Consider the last identity in the lemma. Using the previous identity, it suffices to prove that $(\tau_{|V})_*(c_{|V})=\tau_*(c)_{|\widetilde V}$ for any class $c\in H^*_V(W,\C)$. We can assume that $c$ is represented by a smooth real manifold $Y$ which intersects $V$ transversally. We can also reduce $W$ and $\widetilde W$ in order to assume that there is a projection $\widetilde\Pi:\widetilde W\to \widetilde V$ which defines a smooth fibration whose fibers are diffeomorphic to a ball. We deduce from the discussion before the lemma that 
$$\tau_*(c)\smallsmile \{\widetilde V\}=\tau_*(c)\smallsmile \tau_*\{V\}=\tau_*(c\smallsmile \{V\})=\tau_*\{Y\cap V\}=\{\tau(Y\cap V)\}.$$
Observe that $\tau_*(c)_{|\widetilde V}$ is the image of the class $\tau_*(c)\smallsmile\{\widetilde V\}$ by the natural morphism $\widetilde\Pi_*:H_c^*(\widetilde W,\C)\to H^*_c(\widetilde V,\C)$. The last identities imply that $\tau_*(c)_{|\widetilde V}$ is equal to 
the class of $\tau(Y\cap V)$ in $H^*_c(\widetilde V,\C)$ which is equal to $(\tau_{|V})_*(c_{|V})$ since $c_{|V}$ is represented by $Y\cap V$.
This completes the proof of the lemma.
\endproof

\section{Currents on projective fibrations} \label{section_fibration}

In this section we discuss positive closed currents on fibrations over a complex manifold with projective spaces as fibers. These currents will appear as a kind of derivative in the normal directions along a submanifold, of a positive closed current on a K\"ahler manifold.  Some notions and results can be extended without difficulty to general fibrations where the fibers are not necessarily projective spaces.

Let $V$ be a  K\"ahler manifold of dimension $l$, not necessarily compact, and let $\omega_V$ be a K\"ahler form on $V$. Let $E$ be a holomorphic vector bundle of rank $r$ over $V$ and denote by $\P(E)$ its projectivization. The complex manifold $\P(E)$ is of dimension $l+r-1$. Denote by $\pi:\P(E)\to V$ the canonical projection. The map $\pi$ defines a regular fibration over $V$ with $\P^{r-1}$ fibers.

Consider a Hermitian metric $\|\cdot \|$ on $E$ and denote by  $\omega_{\P(E)}$ the closed $(1,1)$-form on $\P(E)$ induced by $\ddc\log \|v\|$ with $v\in E$. The restriction of $\omega_{\P(E)}$ to each fiber of $\P(E)$ is the Fubini-Study form on this fiber. So $\omega_{\P(E)}$ is strictly positive in the fiber direction.
It follows that given an open set $V_0\Subset V$ there is a constant $c>0$ large enough such that
$c\pi^*(\omega_V)+\omega_{\P(E)}$ is positive on $\pi^{-1}(V_0)$ and defines a K\"ahler metric there.

\begin{definition}\rm \label{def_h_dim}
Let $S$ be a non-zero positive  $(p,p)$-current on $\P(E)$. Let $V_0$ be an open subset of $V$. We call {\it horizontal dimension} (or {\it h-dimension} for short) of $S$ over $V_0$ the largest integer $j$ such that $S\wedge \pi^*(\omega_V^j)\not=0$ on $\pi^{-1}(V_0)$. If this dimension is 0, we say that $S$ is {\it vertical} over $V_0$. 
{\it The h-dimension of $S$} is its h-dimension over $V$. By convention, if $S=0$ on $\pi^{-1}(V_0)$  then the h-dimension of $S$ over $V_0$ is $\max(l-p,0)$.  
\end{definition}

Note that if $s$ is the h-dimension of $S$ then $\max(l-p,0)\leq s\leq \min(l+r-p-1,l)$ because of a dimension reason, see also Lemma \ref{lemma_h_dim_min} below.  
Note also that the positivity of $S$ and the strict positivity of $\omega_V$ imply that the definition does not depend on the choice of $\omega_V$. In fact, we have the following general result.

\begin{lemma} \label{lemma_h_dim}
Let $p$ and $q$ be fixed positive integers. Then the h-dimension of $S$ over $V_0$ is strictly smaller than $\max(p,q)$ if and only if $S\wedge \pi^*(\theta)=0$ for every continuous (or smooth) $(p,q)$-form $\theta$ on $V_0$. 
\end{lemma}
\proof
Define $j:=\max(p,q)$. 
Observe that  $\omega_V^j$ can be written as a finite combination of $\gamma\wedge\theta$ where $\gamma$ is a smooth $(j-p,j-q)$-form and $\theta$ is a smooth $(p,q)$-form. Therefore, the sufficiency of the condition is clear. Assume now that the h-dimension of $S$ over $V_0$ is strictly smaller than $j$. We prove that $S\wedge \pi^*(\theta)=0$. 

Consider first the case where $p=q=j$. Since the problem is local on $V_0$, we can assume that $\theta$ has compact support in $V_0$. Moreover, we can write it as a finite combination of  positive forms with compact support. So we can assume that $\theta$ is positive and $\theta\leq\omega_V^j$. Therefore, we have $0\leq S\wedge \pi^*(\theta)\leq S\wedge \pi^*(\omega_V^j)=0$. It follows that $S\wedge \pi^*(\theta)=0$.

Consider now the case where $(p,q)=(j,j-r)$ with $1\leq r\leq j$. The remaining case can be treated in the same way. 
Observe that $\theta$ can be written as a finite sum of forms of type $\gamma\wedge\beta$ where $\gamma$ is a continuous $(r,0)$-form and $\beta$ is a positive continuous $(j-r,j-r)$-form. So
we can assume that $\theta=\gamma\wedge\beta$. Consider a test smooth form $\theta'$ of appropriate bidegree with compact support in $\pi^{-1}(V_0)$. We have to check that $\langle S\wedge\pi^*(\gamma\wedge\beta),\theta'\rangle=0$. 

As above, we can assume that $\theta'=\gamma'\wedge \beta'$ with $\gamma'$ of bidegree $(0,r)$ and $\beta'$ positive of appropriate bidegree. From Cauchy-Schwarz's inequality, we have 
$$\big|\langle  S\wedge\pi^*(\gamma\wedge\beta),\theta'\rangle\big|\leq  \big|\langle S\wedge\pi^*(\gamma\wedge\overline\gamma\wedge\beta),\beta'\rangle\big|^{1/2}\big|\langle S\wedge\pi^*(\beta),\gamma'\wedge \overline\gamma'\wedge \beta'\rangle\big|^{1/2}.$$
The first factor in the right hand side vanishes according to the bidegree $(j,j)$  case. The result follows.
Note that the same proof holds for $T$ weakly positive and for $(p,q)=(j,j), (j-1,j)$ or $(j,j-1)$. 
\endproof

The following lemma describes the structure of vertical closed currents.

\begin{lemma} \label{lemma_decom_ver}
Let $S$ be a positive closed $(p,p)$-current on $\P(E)$ as above. Assume that $S$ is vertical over an open set $V_0$. Then there exist a unique positive measure $\mu$ on $V_0$ and for $\mu$ almost every $x$, a positive closed 
$(p,p)$-current $S_x$ on $\P(E)$ supported by $\pi^{-1}(x)$ and cohomologous to a linear subspace there, such that 
$$S=\int S_x d\mu(x) \quad \mbox{on}\quad \pi^{-1}(V_0).$$
Moreover, $\mu$ depends linearly on $S$. 
\end{lemma}
\proof
By Lemma \ref{lemma_h_dim}, for any smooth function $\chi$ on $U$ the current $(\chi\circ\pi)S$ is closed on $\pi^{-1}(V_0)$. Since the problem is local on $V_0$, multiplying $S$ with a function $\chi\circ\pi$ with compact support permits to assume that $S$ has support in $\pi^{-1}(K)$ for some compact subset $K$ of $V_0$. Fix a neighbourhood $V_1\Subset V_0$ of $K$. We have seen that $\pi^{-1}(V_1)$ is a K\"ahler manifold. Fix a K\"ahler form on it.

The set of all positive closed $(p,p)$-currents of mass 1 on $\pi^{-1}(V_1)$ satisfying the above property is a convex compact set. Its extremal elements should be supported by a fiber. It follows from Choquet's representation theorem that there is a positive measure $\mu$ on $V$ such that for $\mu$-almost every $x$ there is a positive closed $(p,p)$-current $S_x$ on $\P(E)$ supported by $\pi^{-1}(x)$ such that
$$S=\int S_x d\mu(x).$$
We can multiply $\mu$ with a positive function $\lambda(x)$ and divide $S_x$ by $\lambda(x)$ in order to have that $S_x$ is cohomologous to a linear subspace in $\pi^{-1}(x)$.  We check now that $\mu$ is unique and  depends linearly on $S$.

Fix a closed form $\Omega$ of bidegree $(l+r-1-p,l+r-1-p)$ on $\P(E)$ such that its restriction  to each fiber of $\pi$ is cohomologous to a linear subspace in this fiber, e.g. a power of $\omega_{\P(E)}$. We have
$$S\wedge\Omega = \int (S_x\wedge \Omega) d\mu(x).$$
The intersection $S_x\wedge \Omega$ defines a measure with algebraic mass 1. It follows that $\mu=\pi_*(S\wedge\Omega)$. So $\mu$ is unique and  depends linearly on $S$.
This completes the proof of the lemma.
\endproof

The last lemma and the following one give the complete description of closed currents with minimal h-dimension, i.e. of h-dimension $\max(l-p,0)$.

\begin{lemma} \label{lemma_h_dim_min}
Let $S$ be a positive closed $(p,p)$-current on $\P(E)$ with $p<l$. Assume that the h-dimension of $S$ over $V_0$ is smaller or equal to $l-p$. Then there is a unique positive closed $(p,p)$-current $S^h$ on $V_0$ such that $S=\pi^*(S^h)$ on $\pi^{-1}(V_0)$. In particular, the h-dimension of $S$ over $V_0$ is equal to $l-p$. 
\end{lemma}
\proof
The lemma is clear for $S=0$. So we can assume that $S\not=0$. The uniqueness of $S^h$ is also clear. We prove now the existence of $S^h$. Since this is a local problem, we can assume that $V_0$ is a small ball and $\pi^{-1}(E)$ can be identified with the product $V_0\times\P^{r-1}$, where $\pi$ is identified with the canonical projection on $V_0$. 

Denote by $x=(x_1,\ldots,x_l)$ the complex coordinates on $V_0$. If $I=(i_1,\ldots,i_m)$ with $i_j\in\{1,\ldots,l\}$, define $dx_I:=dx_{i_1}\wedge \ldots\wedge dx_{i_m}$ and 
$d\overline x_I:=d\overline x_{i_1}\wedge \ldots\wedge d\overline x_{i_m}$.
Lemma \ref{lemma_h_dim} implies that $S$ can be written on $V_0\times\P^{r-1}$ as
$$S=\sum_{|I|=|J|=p} R_{IJ} dx_I\wedge d\overline x_J,$$
where $R_{I,J}$ is a 0-current on $V_0\times\P^{r-1}$. 

Since $S$ is closed, we deduce that $dR_{IJ}$ is a combination of $dx_i$ and $d\overline x_i$, i.e. $R_{IJ}$ is constant along the fibers of $\pi$ (to see this point, we can regularize $S$ using some convolution). Therefore, $R_{IJ}$ is the pull-back of a 0-current on $V_0$. We deduce the existence of a current $S^h$ on $V_0$ such that $S=\pi^*(S^h)$. It is clear that $S^h$ should be a positive closed $(p,p)$-current.
\endproof

\begin{proposition} \label{prop_max_h}
Let $S$ be a positive closed $(p,p)$-current on $\P(E)$ as above and let
 $s$ be the h-dimension of $S$ over an open set $V_0$. Let $\Omega$ be a smooth closed form of bidegree $(l-s+r-1-p,l-s+r-1-p)$ on $\pi^{-1}(V_0)$ whose restriction to each fiber of $\pi$ is cohomologous to a linear subspace in this fiber. Then the current $S^h:=\pi_*(S\wedge\Omega)$ on $V_0$ is positive closed of bidegree $(l-s,l-s)$ with support in $\pi(\supp(S))\cap V_0$. Moreover, it does not depend on the choice of $\Omega$. 
\end{proposition}
\proof
It is clear that $S^h$ is a closed $(l-s,l-s)$-current. 
Observe that there is a form $\Omega$ satisfying the hypothesis, e.g. a power of $\omega_{\P(E)}$.
Since the problem is local, we can assume that $V_0$ is strictly contained in $V$ and therefore $\pi^{-1}(V_0)$ is a K\"ahler manifold. In particular, we obtain a strictly positive form $\Omega$ by taking a linear combination of a power of $\omega_{\P(E)}$ and a power of $\pi^*(\omega_V)$.
If we choose $\Omega$ strictly positive, we obtain a positive current $S^h$ with support in $\pi(\supp(S))$.

It remains to prove that $S^h$ does not depend on the choice of $\Omega$.
By Lemma \ref{lemma_h_dim}, if $\alpha$ is a positive closed $(s,s)$-form on $V$ then $S\wedge \pi^*(\alpha)$ is a vertical positive closed $(p+s,p+s)$-current. It follows from this lemma  that $S\wedge \pi^*(\chi\alpha)$ is a vertical positive closed $(p+s,p+s)$-current for any positive function $\chi$ on $V$. 
By Lemma \ref{lemma_decom_ver}, we can associate to this current a measure $\mu$ which depends linearly on $\chi\alpha$. We have seen in the proof of that lemma
that $\mu$ is equal to $\pi_*(S\wedge\pi^*(\chi\alpha)\wedge\Omega)=S^h\wedge\chi\alpha$ and does not depend on the choice of $\Omega$.
Since any $(s,s)$-form $\beta$ can be written as a finite combination of forms of type $\chi\alpha$, the measure $S^h\wedge \beta$ does not depend on the choice of $\Omega$. We deduce that $S^h$ is independent of the choice of $\Omega$. 
\endproof

\begin{definition} \rm
With the notation as in Proposition \ref{prop_max_h}, we say that $S^h$ is {\it the shadow of $S$ on $V_0$}. {\it The shadow of $S$} is its shadow on $V$. 
\end{definition}

Denote by $-h_{\P(E)}$ the tautological class on $\P(E)$ which is the Chern class of the tautological line bundle $O_{\P(E)}(-1)$ over $\P(E)$. With the notation as in the beginning of the section, $h_{\P(E)}$ is the class of $\omega_{\P(E)}$. 
Recall that the cohomology ring $\oplus H_c^*(\P(E),\C)$ is a free  $\oplus H_c^*(V,\C)$-module 
generated by the classes $1,h_{\P(E)},\ldots h_{\P(E)}^{r-1}$, see e.g. Bott-Tu \cite{BottTu} or Voisin \cite{Voisin}.  
This is a consequence of Leray's theorem and can also be deduced from a similar property for de Rham cohomology without compact support via the Poincar\'e duality.
So if $c$ is a class in 
$H_c^{2p}(\P(E),\C)$, we can write it,  in a unique way as
$$c=\sum_{j=\max(0,l-p)}^{\min(l,l+r-1-p)} \pi^*(\kappa_j(c))\smallsmile h_{\P(E)}^{j-l+p},$$
where $\kappa_j(c)$ is a class in $H^{2l-2j}_c(V,\C)$.  If $c$ is the class of a closed $(p,p)$-current $S$ with compact support on $\P(E)$,  we write $\kappa_j(S):=\kappa_j(c)$.

\begin{definition}\rm
The maximal $j$ such that $\kappa_j(c)\not=0$ is called {\it the horizontal dimension} (or {\it h-dimension} for short) of the class $c$. If $c=0$, by convention,  the h-dimension of $c$ is $\max(l-p,0)$.
\end{definition}

The following lemma shows that the h-dimension of a positive closed current with compact support depends only on its cohomology class.
Recall that the class $c$ is said to be {\it pseudo-effective} if it contains a positive closed current.

\begin{lemma} \label{lemma_current_class}
Let $S$ be a positive closed $(p,p)$-current with compact support in $\P(E)$. Then the h-dimension of $S$ is equal to the h-dimension of $\{S\}$. Moreover, if $S^h$ is the shadow of $S$ and $s$ is the h-dimension of $S$ then $S^h$ belongs to the class $\kappa_s(S)$. In particular, $\kappa_s(S)$ is a pseudo-effective class; if $S\not=0$ and if $\pi(\supp(S))$ does not support a positive closed $(l-j,l-j)$-current for some $j\geq 1$, then $s<j$.
\end{lemma}
\proof
Let $\theta$ be a smooth closed $2j$-form on $V$ with $j>s$. By Lemma \ref{lemma_h_dim}, we have $S\wedge \pi^*(\theta)=0$. It follows from the above uniqueness of the decomposition of $H_c^{2p}(\P(E),\C)$ that 
$\kappa_0(S\wedge \pi^*(\theta))=\kappa_j(S)\smallsmile \{\theta\}$. We deduce that $\kappa_j(S)\smallsmile \{\theta\}=0$ for every $\theta$. Hence, by Poincar\'e's duality, $\kappa_j(S)=0$ and the h-dimension of $\{S\}$ is at most equal to $s$. In order to complete the proof of the lemma, we have to check that $S^h$ belongs to $\kappa_s(S)$. 

For this purpose, it is enough to prove that the measure $S^h\wedge \theta$ belongs to $\kappa_0(S\wedge\pi^*(\theta))$ for any smooth closed $2s$-form $\theta$ on $V$. 
We show this property for every smooth $(q,2s-q)$-form $\theta$ not necessarily closed. If $q\not=s$,  we have $S^h\wedge \theta=0$ for bidegree reasons and $S\wedge \pi^*(\theta)=0$ according to Lemma \ref{lemma_h_dim}. So we can assume that $q=s$.
The form $\theta$ can be written as a combination of positive forms. Therefore, we can suppose that $\theta$ is positive. Replacing $S$ with $S\wedge\pi^*(\theta)$ (which is closed according to Lemma \ref{lemma_h_dim}) allows us to assume that $s=0$, i.e. $S$ is a vertical current, and $\theta$ is the constant function 1. So $S^h$ coincides with the measure $\mu$ given in Lemma \ref{lemma_decom_ver}. Using the decomposition given in that lemma, we reduce the problem to the case where $\mu$ is a Dirac mass. We can then check the property without difficulty. 

Note that when $S$ is only weakly positive and $V$ is compact, by Hodge theory, it is enough to consider $\theta$ of bidegree $(j,j)$ or $(s,s)$. So Lemma \ref{lemma_h_dim} still works for these bi-degrees and gives the same result. 
\endproof

We now introduce the notion of $V$-conic currents. They will be used in order to describe the tangent to a positive closed current on a complex manifold along a submanifold $V$.
Let $E$ be a holomorphic vector bundle of rank $r$ over a K\"ahler manifold $V$ of dimension $l$ as above.
We do not assume that $V$ is compact and
we identify it with the zero section of $E$.

The projectivisation $\P(E\oplus \C)$ of the vector bundle $E\oplus \C$ is a natural compactification of $E$. 
Here $\C$ denotes the trivial line bundle over $V$. 
For simplicity, denote $\overline E:=\P(E\oplus \C)$ and $\pi_0:\overline E\to V$ the canonical projection. 
If $V_0$ is an open subset relatively compact in $V$, as we have seen above, $\pi_0^{-1}(V_0)$ is a K\"ahler manifold. 
The action of the multiplicative group $\C^*$ on $E$ extends naturally to $\overline E$.

\begin{definition}\rm
A positive closed $(p,p)$-current $S$ on $E$ is {\it $V$-conic} if it is invariant under the action of $\C^*$.
\end{definition}

We will see in Proposition \ref{prop_conic_char} below that such a current, extended by 0 on the hypersurface at infinity $H_\infty:=\overline E\setminus E$, is a positive closed current on $\overline E$ that we still denote $S$. Note that any current supported by $V$ is $V$-conic. 

Let $\pi_\infty:\overline E\setminus V\to H_\infty$ be the central projection on the hypersurface at infinity. We can also identify $H_\infty$ with $\P(E)$ and the restriction of $\pi_0$ to $H_\infty$ with $\pi:\P(E)\to V$. 
We have the following characterization of $V$-conic currents.

\begin{proposition} \label{prop_conic_char}
Let $S$ be a $V$-conic positive closed $(p,p)$-current as above. Assume that $\supp(S)\cap V$ is compact. Then, there is a unique positive closed $(p,p)$-current $S_\infty$ on $H_\infty\simeq \P(E)$ and a unique positive closed $(p,p)$-current $S_0$ on $\overline E$ with support in $V$ such that 
$$S=\pi_\infty^*(S_\infty)+S_0.$$ 
In particular, $S$ extends by $0$ to a positive closed current on $\overline E$ that we still denote $S$.
Moreover, the intersection $S\wedge [H_\infty]$ is well-defined and is equal to $S_\infty$. The currents $S_\infty,S_0$ have compact supports and $S_0$ is the restriction of $S$ to $V$. 
\end{proposition}

Note that $\pi_\infty^*(S_\infty)$ is well-defined on $\overline E\setminus V$ since $\pi_\infty$ is a submersion there. The following lemma shows that this current extends by 0 to a positive closed $(p,p)$-current on $\overline E$ and we keep the same notation for the extended current. We already obtain here the uniqueness of $S_\infty$. 
The last assertion in the proposition is also clear. Here, by restriction of $S$ to $V$, we mean the multiplication of $S$ with the characteristic function of $V$. One should distinguish it from the intersection of $S$ with $[V]$. 

\begin{lemma} \label{lemma_pull_back}
Let $R$ be a current of order $0$ with compact support on $H_\infty\simeq \P(E)$. Then the current $\pi_\infty^*(R)$ has finite mass near $V$. We still denote by $\pi_\infty^*(R)$ its extension by $0$ through $V$. The operator 
$R\mapsto \pi_\infty^*(R)$ is continuous. If $R$ is closed then $\pi_\infty^*(R)$ is closed and we have $\pi_\infty^*\{R\}=\{\pi_\infty^*(R)\}$. 
\end{lemma}
\proof
Let $\sigma_E:\widehat{\overline E}\to \overline E$ be the blow-up of $\overline E$ along $V$. 
In order to simplify the notation,  we identify $\sigma_E^{-1}(H_\infty)$ with $H_\infty$.
The map $\pi_\infty$ lifts naturally to a {\bf holomorphic} map $\widehat\pi_\infty:\widehat {\overline E}\to H_\infty$ which defines a regular fibration with $\P^1$ fibers over $H_\infty$. So the current $\widehat\pi_\infty^*(R)$ is well-defined and of order 0 and has no mass on $\sigma_E^{-1}(V)$. Therefore, $(\sigma_E)_*\widehat\pi_\infty^*(R)$ is a current of order 0 having no mass on $V$. It is equal to $\pi_\infty^*(R)$ outside $V$ and depends continuously on $R$. 
The first and second assertions in the lemma follow.

For the last assertion, assume that $R$ is closed. Clearly, $\pi_\infty^*(R)$ is closed. It remains to check the identity in the lemma. Since $\pi_\infty^*(R)$ depends continuously on $R$, by de Rham's regularization theorem, it is enough to consider the case where $R$ is a smooth form. Recall that the operator 
$\pi_\infty^*:H^*_c(H_\infty,\C)\to H^*_c(\P(E),\C)$ is defined by $\pi_\infty^*\{R\}:=\{(\sigma_E)_*\widehat\pi_\infty^*(R)\}$ for $R$ smooth and closed. Note that even when $R$ is smooth $(\sigma_E)_*\widehat\pi_\infty^*(R)$ should be considered as a current.
The key point here is that de Rham cohomology groups can be defined using smooth forms or currents.
The definition of the action of $\pi_\infty^*$ on cohomology is in fact valid for more general meromorphic maps. 
The last identity in the lemma is clear for smooth forms $R$.
\endproof

We also need the following lemma. Assume that $V$ is a submanifold of a K\"ahler manifold $X$ of dimension $k$. We use the notations introduced in Section \ref{section_preliminar}. Recall that 
 $\sigma:\widehat X\to X$ is the blow-up of $X$ along $V$ and $\widehat V:=\sigma^{-1}(V)$.
 
\begin{lemma} \label{lemma_strict_transform}
Let $T$ be a positive closed $(p,p)$-current on $X$ with support in a fixed open set $W_1$ of $X$ such that $W_1\cap V\Subset V$.  Then, for every open sets $U\Subset U'\Subset X$ containing $W_1\cap V$, the mass of 
$(\sigma_{|\widehat X\setminus \widehat V})^*(T)$ on $\sigma^{-1}(U)\setminus \widehat V$ is bounded by $c\|T\|_{U'}$ for some constant $c>0$ independent of $T$. In particular, $(\sigma_{|\widehat X\setminus \widehat V})^*(T)$ extends by $0$ to a positive closed $(p,p)$-current on $\widehat X$ that we denote by $\sigma^\diamond(T)$.
\end{lemma}
\proof
The second assertion is a consequence of the first one and of an extension theorem by Skoda \cite{Skoda}. Let $\widehat\omega$ be a K\"ahler form on $\widehat W_1:=\sigma^{-1}(W_1)$, see also Section \ref{section_preliminar}.
Observe that $R:=\sigma_*(\widehat\omega)$ is a positive closed $(1,1)$-current which is smooth outside $V$ and has no mass on $V$.  The mass of $(\sigma_{|\widehat X\setminus \widehat V})^*(T)$ on $\sigma^{-1}(U)\setminus \widehat V$ is equal to the mass of the measure $T\wedge R^{k-p}$ on $U\setminus V$. 
We have to bound  the last quantity.

Let $\varphi$ be the quasi-psh function on $W_1$ introduced in Section \ref{section_preliminar} such that 
$\ddc\varphi+c_0\omega\geq R$. This function is smooth outside $V$.
Define for $M>0$ large enough $\varphi_M:=\max(\varphi,-M)$ and $\omega_M:=c_0\omega+\ddc\varphi_M$. It is not difficult to see that $\omega_M$ is a positive closed current  which is larger than $R$ on the open set $\{\varphi> -M\}$. When $M\to\infty$, this open set increases to $W_1\setminus V$.
Therefore, it is sufficient to bound the mass of $T\wedge \omega_M^{k-p}$ on $U$ by a quantity which is independent of $M$ for $M$ large enough.

Observe that the positive measure $T\wedge \omega_M^{k-p}$ is well-defined because $\omega_M$ has continuous local potentials. Moreover, for $M$ large enough, $\varphi_M=\varphi$ on $W_1\cap U'\setminus U$. It follows from Stokes' formula that the mass of $T\wedge \omega_M^{k-p}$ on $U$ does not depend on $M$. 
Fix an $M$ large enough. The classical Chern-Levine-Nirenberg's inequality implies that this mass is bounded by a constant times $\|T\|_{U'}$, see Chern-Levine-Nirenberg \cite{ChernLevineNirenberg} and  Demailly \cite{Demailly3}. The lemma follows.
\endproof

Note that the result can be generalized for maps between manifolds of different dimensions, see also \cite{DinhSibony8}. In this paper, we only need the version stated above.

\begin{definition} \rm
With the notations as in Lemma \ref{lemma_strict_transform}, we call $\sigma^\diamond(T)$ {\it the strict transform} of $T$ by $\sigma$.
\end{definition}

In general, $\sigma^*\{T\}$ is not equal to $\{\sigma^\diamond(T)\}$. The missing part is described in the following lemma.

\begin{lemma} \label{lemma_exceptional_class}
With the notations of Lemma \ref{lemma_strict_transform}, there is a class $e(T)$ in $H^{2p-2}_c(\widehat V,\C)$ such that for any neighbourhood $\widehat W$ of $\widehat V$ the class $\sigma^*\{T\}-\{\sigma^\diamond(T)\}$ in $H^{2p}_{\widehat V}(\widehat W,\C)$ is equal to the canonical image of $e(T)$ in this cohomology group. 
\end{lemma}
\proof
Choose a neighbourhood $\widehat W'$ of $\widehat V$ which is $\widehat V$-contractile, i.e. there is a smooth projection $\Pi:\widehat W'\to \widehat V$ which defines a fibration with connected and simply connected fibers. By de Rham's regularization theorem, there is a current $T'$ with support in $W_1$ smooth near $\widehat V$ and equal to $T$ outside a compact set in $W':=\sigma(\widehat W')$  such that the class of $T-T'$ in $H^{2p}_c(W',\C)$ vanishes. Define $e(T)$ as the class of the current $\Pi_*(\sigma^*(T')-\sigma^\diamond(T))$. It is clear that the property in the lemma is true for $\widehat W'$. Observe also that the $e(T)$ does not depend on the choice of $T'$.

Consider now an arbitrary open set $\widehat W$ as in the lemma. Choose a neighbourhood $\widehat W''\subset \widehat W\cap \widehat W'$ of $\widehat V$ such that $\Pi$ restricted to $\widehat W''$ defines a fibration with connected and simply connected fibers.  Since $e(T)$ does not depend on the above choice of $T'$, we can choose a $T'$ such that 
$T-T'$ is supported by $W'':=\sigma(\widehat W'')$ and its class in $H^*_c(W'',\C)$ vanishes. As above, we see that the property in the lemma holds for $\widehat W$.
\endproof

\noindent
{\bf Proof of Proposition \ref{prop_conic_char}.}
It is well-know that we can decompose $S$ in a unique way into a sum of two positive closed $(p,p)$-currents
$S=S'+S_0$ with $S_0$ supported by $V$ and $S'$ without mass on $V$, see e.g. Demailly \cite{Demailly3} and Skoda \cite{Skoda}. For simplicity, we replace $S$ with $S'$ in order to assume that $S_0=0$ and $S$ has no mass on $V$.

Consider first the case where $V$ is a hypersurface. So $\pi_\infty$ extends to a holomorphic map on $\overline E$ and defines a regular fibration over $H_\infty$ with $\P^1$ fibers. Locally, we can identify this fibration with the product $B\times \P^1$ where $B$ is a ball in $\C^{k-1}$. The hypersurface $H_\infty$ is identified with $B\times\{\infty\}$ and the map $\pi_\infty$ is identified with the canonical projection on $B$. The hypersurface $V$ is identified with $B\times \{0\}$. 

Since $S$ is invariant under the action of $\C^*$, we can write on $B\times \C^*$ using the natural coordinates $(z,t)$ 
$$S=S_1(z)\wedge {idt\wedge d\overline t\over t^2} +S_2(z)\wedge {dt\over t} +\overline {S_2(z)}\wedge {d\overline t\over \overline t} + S_3(z),$$
where the $S_i$ are currents of order 0 and of the appropriate bidegree which do not depend on $t$.
Since $S$ has finite mass near $B\times \{0\}$, the first term vanishes. Then the positivity of $S$ implies that the next two terms vanish. This implies the proposition for the hypersurface case with $S_\infty$ such that $S_3=\pi_\infty^*(S_\infty)$.

Consider now the general case. We use the notation introduced in Lemma \ref{lemma_pull_back}. 
Let $\widehat S$ be the strict transform of $S$ by $\sigma_E$ (we use here the hypothesis on the support of $S$). 
The action of $\C^*$ on $E$ extends to $\overline E$ and can be lifted to $\widehat{\overline E}$. 
The current $\widehat S$ is still invariant under this action. We can apply the hypersurface case considered above to the current $\widehat S$ and to the exceptional hypersurface $\sigma_E^{-1}(V)$. 
As in Lemma \ref{lemma_pull_back}, for simplicity, we identify $\sigma_E^{-1}(H_\infty)$ with $H_\infty$. 
So we can write $\widehat S=\widehat\pi_\infty^*(S_\infty)$ with a positive closed $(p,p)$-current $S_\infty$ on $H_\infty$.  Since these currents have no mass on $V$, we deduce that $S=\pi_\infty^*(S_\infty)$.  This completes the proof of the proposition.
\hfill $\square$

\medskip

Let $S$ be a $V$-conic current as above with compact support in $\overline E$. Let $-h_{\overline E}$ denote the tautological class of $\overline E=\P(E\otimes \C)$. 
By Leray's theorem, we can decompose the class of $S$ as
$$\{S\}=\sum_{j=\max(0,l-p)}^{\min(l,l+r-p)} \pi_0^*(\kappa_j(S))\smallsmile h_{\overline E}^{j-l+p}$$
where $\kappa_j(S)$ is a class in $H^{2l-2j}_c(V,\C)$ with $\kappa_j(S)=0$ when $j$ does not satisfies the inequalities 
$\max(0,l-p)\leq j\leq \min(l,l+r-p)$. 
We have the following lemma.

\begin{lemma} \label{lemma_compare_infty}
Let $S,S_\infty$ and $S_0$ be as in Proposition \ref{prop_conic_char}. Then 
$$\kappa_{l+r-p}(S)=\{S_0\} \quad \mbox{and} \quad \kappa_j(S)=\kappa_j(S_\infty) \quad \mbox{for} \quad j<l+r-p,$$ 
where $\{S_0\}$ is the class of $S_0$ in $H^{2p-2r}_c(V,\C)$. 
In particular, if $s$ is the h-dimension of $S_\infty$, then $\kappa_j(S)=0$ for $j>s$ except possibly for $j=l+r-p$ and $\kappa_s(S)$ contains the shadow of $S_\infty$. 
\end{lemma}
\proof
Observe that the second assertion is a consequence of the first one and of Lemma \ref{lemma_current_class}.
Recall that for simplicity we identify $H_\infty$ with $\P(E)$ and $\overline E$ with $\P(E\oplus \C)$. The map $\pi_\infty$ is induced by the canonical projection from $E\oplus \C$ to $E$. The pull-back of a Hermitian metric on $E$ gives a singular Hermitian metric on $E\oplus \C$. We deduce that $\pi_\infty^*(h_{\P(E)})=h_{\overline E}$. 
Using the blow-up as in Lemma \ref{lemma_pull_back}, we obtain easily that 
$$\pi_\infty^*(h_{\P(E)}^m)=h_{\overline E}^m \quad \mbox{for} \quad m<r \qquad \mbox{and} \qquad h_{\overline E}^r=[\pi_\infty^*(h_{\P(E)})]^r=\{V\}.$$

Therefore, using that $\pi_0$ is identified with $\pi\circ\pi_\infty$, we get
\begin{eqnarray*}
\{\pi_\infty^*(S_\infty)\} & = & \sum_{j=\max(0,l-p)}^{\min(l,l+r-1-p)} \pi_\infty^*\pi^*(\kappa_j(S_\infty))\smallsmile \pi_\infty^*(h_{\P(E)}^{j-l+p})\\
& = & \sum_{j=\max(0,l-p)}^{\min(l,l+r-1-p)} \pi_0^*(\kappa_j(S_\infty))\smallsmile h_{\overline E}^{j-l+p}
\end{eqnarray*}
and the class of $S_0$ in $H^{2p-2r}_V(\overline E,\C)$ is equal to
$$\pi_0^*\{S_0\}\smallsmile \{V\} = \pi_0^*\{S_0\}\smallsmile h_{\overline E}^r.$$
Then, the lemma follows from the identity 
$S=\pi_\infty^*(S_\infty)+S_0$ and the uniqueness of the above decompositions. 
\endproof

In the following lemma, we can use any fixed Hermitian metric on $\overline E$. 

\begin{lemma} \label{lemma_conic_mass}
Let $K$ be a fixed compact subset of $V$ and let $U$ be a fixed neighbourhood of $K$ in $E$. If $S$ is a $V$-conic $(p,p)$-current with support in $\pi_0^{-1}(K)$, then $\|S\|\leq c\|S\|_U$ for some constant $c>0$ independent of $S$. 
\end{lemma}
\proof
If the lemma were wrong, there would be a sequence of $V$-conic $(p,p)$-currents $(S_n)$ supported by $\pi_0^{-1}(K)$ such that $\|S_n\|\geq n\|S_n\|_U$. We can divide each $S_n$ by its mass in order to assume that $\|S_n\|=1$. Extracting a subsequence allows to assume that $S_n$ converges to a $V$-conic current $S$ of mass 1 which vanishes on $U$. Since this current is $V$-conic, it vanishes on $E$. This is a contradiction. 
\endproof

We come back to the case where $V$ is a submanifold of dimension $l$ of a K\"ahler manifold $X$ of dimension $k$. Let $\sigma: \widehat X\to X$ and $\widehat V:=\sigma^{-1}(V)$ be as above.  Denote by $E$ the normal vector bundle to $V$ in $X$.
Then the exceptional hypersurface $\widehat V$ of $\widehat X$ is canonically  identified with $\P(E)$. So we can identify the restriction of $\sigma$ to $\widehat V$ with $\pi:\P(E)\to V$. We will need the following lemma.

\begin{lemma} \label{lemma_negative_intersection}
Let $S$ be a positive closed $(p,p)$-current on $\widehat X$ with compact support in $\widehat V=\P(E)$ and with $p\geq 1$. Let $s$ be the h-dimension of $S$. 
Assume that $s$ is strictly smaller than the complex dimension $k-p$ of $S$. Let $\{S\}'$ denote the class of $S$ in $H^{2p}_{\widehat V}(\widehat X,\C)$. 
Then $\kappa_j(\{S\}'_{|\widehat V})=0$ if $j>s$ and  $-\kappa_s(\{S\}'_{|\widehat V})$ contains the shadow of $S$ on $V$. 
In particular, the class $-\kappa_s(\{S\}'_{|\widehat V})$ is  pseudo-effective.
\end{lemma}
\proof
Using a diffeomorphism from a neighbourhood of $\widehat V$ in $\widehat X$ to a neighbourhood of $\widehat V$ in $\widehat{\overline E}$ which is identity on $\widehat V$, we can reduce the problem to the case where $X=\overline E$ and $\widehat X=\widehat{\overline E}$. Let $\widehat\pi_0:\widehat{\overline E}\to \widehat V$ be the canonical projection. It defines a fibration with $\P^1$ fibers over $\widehat V$. 

We can identify $S$ with the intersection of $\widehat\pi_0^*(S)$ with $[\widehat V]$. 
Let $\{S\}$ denote the class of $S$ in $H^{2p-2}_c(\widehat V,\C)$. 
We have 
$$\{S\}'_{|\widehat V}=(\widehat\pi_0^*\{S\}\smallsmile [\widehat V])_{|\widehat V}=\{S\}\smallsmile [\widehat V]_{|\widehat V}=-\{S\}\smallsmile h_{\P(E)}.$$
Finally, since $s$ is strictly smaller than the complex dimension of $S$, we deduce from the definition of $\kappa_j(\cdot)$ that $-\kappa_j(\{S\}'_{|\widehat V})=\kappa_j(\{S\})$. The lemma follows.
\endproof

\section{Tangent cones for positive closed currents} \label{section_tangent}

In this section, we introduce  the tangent cones, along a submanifold, of a positive closed current on a K\"ahler manifold.  We refer the reader to Siu \cite{Siu} for the case where the submanifold is just a point, i.e. the Lelong number case.

Let $X$ be a K\"ahler manifold of dimension $k$. Let $V$ be a submanifold of dimension $l$.
Let $T$ be a positive closed $(p,p)$-current on $X$ such that $\supp(T)\cap V$ is compact. The later condition is satisfied when $V$ or $X$ is already compact. We want to define tangent currents to $T$ along $V$. They are $V$-conic currents on $\overline E$ where $E:=N_{V|X}$ is the normal vector bundle to $V$ in $X$.
We need a special class of homeomorphisms from neighbourhoods of $V$ in $X$ to neighbourhoods of $V$ in $E$ which are in some sense close to holomorphic maps near $V$.

Consider a point $a$ in $V$. If $\Tan_aX$ and $\Tan_aV$ denote respectively  the tangent spaces of $X$ and of $V$ at  $a$, the fiber $E_a$  of $E$ over $a$ is canonically identified with the quotient space $\Tan_aX/\Tan_a V$. Let $x=(x',x'')$ with $x'=(x_1,\ldots,x_l)$ and $x''=(x_{l+1},\ldots,x_k)$ be a local holomorphic coordinate system that identifies a chart of $X$ to the polydisc $2\D^k$ in $\C^k$ such that $V$ is defined by the equation $x''=0$ in this polydisc. In these local coordinates, the bundle $E$ is canonically identified over $V\cap 2\D^k$ with the trivial bundle $(V\cap 2\D^k)\times \C^{k-l}$ which is an open subset of $\C^k$ containing $2\D^k$.

Let $V_0$ be an open subset of $V$. Let $\tau$ be a bi-Lipschitz map from a neighbourhood of $V_0$ in $X$ to a neighbourhood of $V_0$ in $E$. We assume that $\tau$ is equal to identity on $V_0$ and is smooth outside $V_0$. 

\begin{definition} \rm \label{def_admissible}
We say that $\tau$ is {\it admissible} (resp. {\it almost-admissible}) if in any local holomorphic coordinate system as above $\tau$ is admissible (resp. almost-admissible) in the sense of Definition \ref{def_admissible_loc} (resp. Definition \ref{def_almost_admissible_loc}). 
\end{definition}

Note that if $\tau$ is smooth and admissible its differential at any point of $V$ is $\C$-linear and induces an endomorphism of $E$ which is equal to identity. 
If $V_0$ is small enough we can find $\tau$ admissible and holomorphic. In general, we have the following lemma, see also Lemma \ref{lemma_lift_admissible}.

\begin{lemma} \label{lemma_admissible}
There is a smooth admissible map for $V_0=V$. 
\end{lemma}
\proof
Consider a Hermitian metric on $X$. It induces a Hermitian metric on the tangent bundle of $X$. Denote by $F$ the restriction of this tangent bundle to $V$. The tangent bundle of $V$ is identified with a vector sub-bundle $F'$ of $F$ and $E$ is identified with $F/F'$. Let $F''$ denote the orthogonal complement of $F'$ in $F$. This is a vector bundle over $V$ with complex fibers but in general it is {\bf not a  holomorphic} vector bundle. 

The canonical projection $\tau_1:F''\to E$ is a smooth isomorphism between vector bundles and is $\C$-linear on each fiber. Let $\tau_2:F''\to X$ be the map induced by the exponential maps at the points of $V$. This map defines a smooth diffeomorphism between a neighbourhood of $V$ in $F''$ and a neighbourhood of $V$ in $X$. It is identity on $V$ and its differential at each point of $V$ is identity. Define $\tau:=\tau_1\circ \tau_2^{-1}$ on a small neighbourhood of $V$ in $X$. This is a diffeomorphism between this neighbourhood and its image which is a neighbourhood of $V$ in $E$. 

In local coordinates $x=(x',x'')$ as above, we can find smooth matrix-functions $a(x')$, such that the fiber of $F''$ over a point $(x',0)$ is the set of points $(x'+x''a(x'),x'')$ with $x''\in \C^{k-l}$. In a small neighbourhood of $V$, these affine spaces are pairwise disjoint. The map $\tau_1$ sends $(x'+x''a(x'),x'')$ to $(x',x'')$. 
Hence, it sends  $(x',x'')$ to $(x'-x''a(x')+O(\|x''\|^2),x'')$. On the other hand,  
the map $\tau_2^{-1}$ is smooth and tangent to identity at each point of $V$. So we have $\tau_2^{-1}(x)=x+O(\|x''\|^2)$ and 
$d\tau_2^{-1}(x)=dx+O^{**}(\|x''\|^2)$.
So $\tau$ satisfies Definition \ref{def_admissible_loc}. 
\endproof

In what follows, we often use the blow-up $\sigma:\widehat X\to X$ of $X$ along $V$ and the blow-up 
$\sigma_{E}:\widehat {\overline E}\to \overline E$  of $\overline E$ along $V$.  
We will show that some admissible maps on $X$ can be lifted to almost-admissible maps on $\widehat X$. However, in general, we loose the smoothness of these maps and they are only bi-Lipschitz. This is the motivation for Definition \ref{def_admissible_loc}. 

Observe that $\sigma_{E}^{-1}(V)$ can be canonically identified with $\P(E)$. So we also identify it with $\widehat V$. The restriction of $\sigma_{E}$ to this hypersurface is identified with the restriction of $\sigma$ to $\widehat V$ and with $\pi:\P(E)\to V$. For simplicity, we identify $\sigma_{E}^{-1}(H_\infty)$ with $H_\infty$. The projections $\pi_0:\overline E\to V$ and $\pi_\infty:\overline E\setminus V\to H_\infty$ lift to projections $\widehat\pi_0:\widehat {\overline E}\to\widehat V$ and $\widehat \pi_\infty:\widehat {\overline E}\setminus\widehat V\to H_\infty$ in a canonical way. We have $\sigma_E\circ \widehat\pi_0=\pi_0\circ \sigma_E$ and $\widehat\pi_\infty=\pi_\infty\circ \sigma_E$. 
The map $\widehat\pi_\infty$ extends holomorphically to $\widehat{\overline E}$. Both $\widehat\pi_0$ and $\widehat\pi_\infty$ define fibrations with $\P^1$ as fibers. 
Finally, $\widehat E:=\widehat{\overline E}\setminus H_\infty$ can be identified with the normal line bundle to $\widehat V$ in $\widehat X$ and $\widehat {\overline E}$ is its natural compactification. 

\begin{lemma} \label{lemma_lift_admissible}
Let $\tau$ be the smooth admissible map  constructed in Lemma \ref{lemma_admissible}. Then there is a unique almost-admissible map $\widehat \tau$, from a neighbourhood of $\widehat V_0:=\sigma^{-1}(V_0)$ in $\widehat X$ to a neighbourhood of $\widehat V_0$ in $\widehat E$ 
such that $\sigma_{E}\circ\widehat\tau=\tau\circ\sigma$. 
\end{lemma}
\proof
Since $\sigma$ and $\sigma_E$ are  biholomorphic maps outside $\widehat V$, we necessarily have $\widehat\tau=\sigma_{E}^{-1}\circ\tau\circ\sigma$ outside $\widehat V$. By continuity, the map $\widehat\tau$ is unique if it exists. We will describe $\widehat\tau$ outside $\widehat V$ using local coordinates and we will see that it extends to an almost-admissible map.

Let $x=(x',x'')$ be as above where we identify a chart of $X$ with $2\D^k$. The restriction of $V$ to $2\D^k$ is given by the equation $x''=0$. The vector bundle $E$ is identified over $V\cap 2\D^k$ with $2\D^l\times \C^{k-l}$. The map $\tau$ is described as in Definition \ref{def_admissible_loc}. In these coordinates, we identify $X$ with $E$ and $\sigma$ with $\sigma_{E}$ over $2\D^k$. 
Consider the chart $\widehat D$ of $\sigma^{-1}(2\D^k)$ as in Section \ref{section_preliminar} with coordinates $z=(z_1,\ldots, z_k)$ such that $|z_j|< 2$ and
$$\sigma(z)=(z_1,\ldots,z_l,z_{l+1}z_k,\ldots, z_{k-1}z_k,z_k).$$

In these coordinates, $\widehat V$ is given by $z_k=0$. We have
$$\sigma^{-1}_E(x)=\sigma^{-1}(x)=(x_1,\ldots,x_l,x_{l+1}x_k^{-1},\ldots,x_{k-1}x_k^{-1},x_k).$$
Define $z':=(z_1,\ldots,z_l)$ and $z^\#:=(z_{l+1},\ldots,z_{k-1})$.
Using the local description of $\tau$ in the proof of Lemma \ref{lemma_admissible}, we can write the map $\widehat\tau:=\sigma_{E}^{-1}\circ \tau\circ\sigma$ on $\widehat D\setminus\widehat V$  in coordinates $z$ as
$$\widehat\tau(z)=\big(z'+ O^*(|z_k|),z^\#+O^*(|z_k|),z_k+O^*(|z_k|^2)\big).$$
We see that $\widehat\tau$ extends continuously to a map on $\sigma^{-1}(\D^k)$ which is identity on $\widehat V$. 

In the last identity, the function $O^*(|z_k|^2)$ is smooth and $O^*(|z_k|)$ is the product of $z_k^{-1}$ with a smooth $O(|z_k|^2)$ function. Since smooth $O(|z_k|^2)$ functions can be written as a combination of $z_k^2,z_k\overline z_k$ and $\overline z_k^2$ with smooth coefficients, we see that
$$d\widehat\tau(z)=\big( dz'+O^{**}(|z_k|),dz^\#+O^{**}(|z_k|), dz_k + O^{**}(|z_k|^2)\big).$$
Hence, $\widehat\tau$ is almost-admissible for $\widehat V$ in $\widehat X$.
\endproof

Denote by $A_\lambda$ the automorphism of $\overline E$ induced by the multiplication by $\lambda\in\C^*$. 
Let $\tau$ be an {\bf almost-admissible} map as in Definition \ref{def_admissible}. 
Fix an open subset $W_1$ of $X$ such that $W_1\cap V$ is non-empty and relatively compact in $X$. Consider a positive closed $(p,p)$-current on $X$ with support in $W_1$. 
We can decompose $T$ as $T'+T_0$ where $T',T_0$ are positive closed currents, $T'$ has no mass on $V$ and $T_0$ is the restriction of $T$ to $V$.

By Lemma \ref{lemma_bi_lipschitz}, we can define a closed $2p$-current $\tau_*(T')$ of order 0 on a neighbourhood of $V_0$ in $E$ with no mass on $V$. Since $\tau$ is identity on $V$, we define $\tau_*(T)=\tau_*(T')+T_0$.   Define $T_\lambda:=(A_\lambda)_*(\tau_*(T))$. This is a closed $2p$-current whose domain of definition converges to an open set containing $\pi_0^{-1}(V_0)\setminus H_\infty$.

\begin{proposition} \label{prop_existence_tangent}
Let $U$ be an open set relatively compact in $\pi_0^{-1}(V_0)\setminus H_\infty$. Then for $\lambda$ large enough, the current $T_\lambda$ is defined on $U$ and its mass on $U$ is bounded by $c\|T\|$ for some constant $c$ independent of $\lambda$ and $T$. Moreover, 
if $(\lambda_n)$ is a sequence converging to infinity such that $T_{\lambda_n}$ converges to a current $S$ on $\pi_0^{-1}(V_0)\setminus H_\infty$, then $S$ is a positive closed $(p,p)$-current independent of the choice of $\tau$.
\end{proposition}
\proof
We only consider $\lambda$ large enough. So
the current $T_\lambda$ is closed and is defined on $U$.
Since the problem is local with respect to $V$, we can assume that $V_0$ is equal to $\D^k\cap V$, where $\D^k$ is identified with a chart of $X$ and a chart of $E$ with local coordinates $x$ as above. 

Let $R$ be a smooth $(2k-2p)$-form with compact support in $U$ and with coefficients bounded by 1. 
We have $\langle T_\lambda,R\rangle=\langle T,\tau^*(A_\lambda)^*(R)\rangle$. By Remark \ref{rk_almost_admissible}, the last integral is bounded by a constant times $\|T\|$. It follows that the mass of $T_\lambda$ on $U$ is bounded by a constant times $\|T\|$. In particular, for any sequence $(\lambda_n)$ converging to infinity, we can extract a subsequence $(\lambda_{n_i})$ such that $T_{\lambda_{n_i}}$ converges to a closed current on $\pi_0^{-1}(V_0)\setminus H_\infty$.

Let $R'$ denote the component of bidegree $(k-p,k-p)$ of $R$. Define $R_\lambda:=(A_\lambda)^*(R')$. By Remark \ref{rk_almost_admissible}, $\langle T_\lambda,R\rangle-\langle T,R_\lambda\rangle$ tends to 0 as $\lambda$ tends to infinity. 
Since $R_\lambda$ does not depend on the choice of $\tau$, we deduce that $S$ does not depend on the choice of $\tau$.
If $R$ is of bidegree $(q_1,q_2)$ with $(q_1,q_2)\not=(k-p,k-p)$ then $R'=0$. It follows that $\langle S,R\rangle=0$. Hence, $S$ is a current of bidegree $(p,p)$. When $R$ is a weakly positive $(k-p,k-p)$-form, $R_\lambda$ is also weakly positive and hence $\langle T,R_\lambda\rangle$ is positive. We deduce that $\langle S,R\rangle$ is positive. Therefore, $S$ is a positive current.
\endproof

Let $\tau$ be a global {\bf almost-admissible} map as in Definition \ref{def_admissible} for $V_0=V$. 
We define as above $T_\lambda:=(A_\lambda)_*\tau_*(T)$. By Proposition \ref{prop_existence_tangent}, the following notion of tangent current does not depend on the choice of $\tau$.

\begin{definition} \rm
A current $S$ on $E$ is said to be a {\it tangent current} to $T$ along $V$ if there is a sequence $\lambda_n\to\infty$ such that $S=\lim T_{\lambda_n}$. We also say that $S$ is the tangent current associated with the sequence $(\lambda_n)$. 
\end{definition}

Observe that if $\theta$ is a smooth positive closed $(q,q)$-form on $X$ and if $S$ is the tangent current to $T$ along $V$ associated with a sequence $(\lambda_n)$ then $S\wedge\pi_0^*(\theta_{|V})$ is the tangent current to $T\wedge\theta$ along $V$ associated with $(\lambda_n)$.

\begin{theorem} \label{th_tangent_conic}
Let $X$ be a K\"ahler manifold and let $V$ be a submanifold of $X$. Denote by $E$ the normal vector bundle to $V$ in $X$,  $\overline E$ its natural compactification and $\pi_0:\overline E\to V$ the canonical projection. Let $W_1$ be a fixed open subset of $X$ such that $W_1\cap V$ is relatively compact in $X$. 
If $T$ is a positive closed $(p,p)$-current on $X$ with support in $W_1$, then its
tangent currents along $V$ are $V$-conic supported by $\pi_0^{-1}(\supp(T)\cap V)$ and of mass bounded by $c\|T\|$ for some constant $c$ independent of $T$. Moreover, these tangent currents belong to the same cohomology class in $H^{2p}_c(\overline E,\R)$ and their restrictions to $V$ are equal to the restriction of $T$ to $V$.
\end{theorem}

We need the following lemma where $\tau$ is {\bf smooth admissible} and $\widehat\tau$ is given by Lemma \ref{lemma_lift_admissible}.

\begin{lemma} \label{lemma_tangent_blowup}
Let $S$ be the tangent current to $T$ along $V$ associated with a sequence $(\lambda_n)$. 
Then the restriction of $S$ to $V$ is equal to the restriction of $T$ to $V$. 
Let $\widehat T$ be the strict transform of $T$ by $\sigma:\widehat X\to X$. Then $\widehat T$ admits a tangent current $\widehat S$ along $\widehat V$ associated with the same sequence $(\lambda_n)$. Moreover, $\widehat S$ is the strict transform of $S$ by $\sigma_{E}:\widehat{\overline E}\to \overline E$.  
\end{lemma}
\proof
If $T$ is supported by $V$ then $\widehat T=0$ and $S=T$. The lemma is clear. So we can assume that $T$ has no mass on $V$. 
We can replace $(\lambda_n)$ by a subsequence in order to assume that $\widehat T$ admits a tangent current  $\widehat S$ along $\widehat V$ associated with $(\lambda_n)$. We have to check that it is the strict transform of $S$ and that $(\sigma_{E})_*(\widehat S)=S$. The last equality implies that $S$ has no mass on $V$.

Denote by $\widehat A_\lambda$ the map on $\widehat {\overline E}$ induced by the multiplication by $\lambda$. It is the natural lift of $A_\lambda$ to $\widehat{\overline E}$. Define also $\widehat T_\lambda:=(\widehat A_\lambda)_*\widehat\tau_*(\widehat T)$. 
Let $R$ be a test smooth $(2k-2p)$-form with compact support in $E$. Define $\widehat R:=\sigma_{E}^*(R)$. It is not difficult to see that 
$\langle \widehat T_\lambda,\widehat R\rangle =\langle T_\lambda,R\rangle$ because $T_\lambda$ and $\widehat T_\lambda$ have no mass on $V$ and $\widehat V$ respectively.
It follows that $(\sigma_{E})_*(\widehat S)=S$. Since $\sigma_{E}$ is injective outside $\widehat V$, it remains to check that $\widehat S$ has no mass on $\widehat V$. 

In order to simplify the notation, we consider the case where $V$ is a hypersurface and $T$ has no mass on $V$. We have to check that $S$ has no mass on $V$. The result we obtain when applied to $\widehat X,\widehat V,\widehat T$ and $\widehat S$, gives the lemma. Multiplying $T$ with a strictly positive closed form allows us to reduce the problem to the case where $T$ is of bidegree $(k-1,k-1)$, see the observation before Theorem \ref{th_tangent_conic}. 

We use local coordinates $x=(x',x_k)$ with $x'=(x_1,\ldots,x_{k-1})$ on a chart $\D^k$ as above. 
Let $\gamma$ denote the restriction of $\ddc\|x'\|^2$ to 
 $\D^k$. 
The mass of $S$ on $V\cap \D^k$ is bounded by a constant times
$\langle S, \gamma\rangle$. Arguing as in Proposition \ref{prop_existence_tangent}, we see that the last integral is bounded by
$$\limsup_{\lambda\to\infty} \langle T, (A_\lambda)^*(\gamma)\rangle  = 
\limsup_{\lambda\to\infty} \langle T, \ddc\|x'\|^2\rangle_{\D^l\times\lambda^{-1}\D^{k-l}}=0$$
since $T$ has no mass on $V$. It follows that $S$ has no mass on $V$. The proof of the lemma is now complete.
\endproof

\noindent
{\bf End of the proof of Theorem \ref{th_tangent_conic}.} 
The theorem is clear when $T$ is supported by $V$. So  we can assume that $T$ has no mass on $V$. The last assertion is already obtained in Lemma \ref{lemma_tangent_blowup}. The mass estimate for tangent currents is a consequence of Proposition \ref{prop_existence_tangent} and Lemma \ref{lemma_conic_mass}. 
The assertion on the supports of the tangent currents is also clear.
We prove now that the tangent currents are $V$-conic and that they have the same cohomology class. By Lemma \ref{lemma_tangent_blowup}, we can assume that $V$ is a hypersurface of $X$.

We use a chart $\D^k$ of $X$ as above with local coordinates $x=(x_1,\ldots,x_k)$ such that $V\cap\D^k$ is given by the equation $x_k=0$. Let $R$ be a smooth $(k-p,k-p)$-form with compact support in $\D^k$. With notations as above, we have seen in Proposition \ref{prop_existence_tangent} that $\langle T_\lambda,R\rangle - \langle T, (A_\lambda)^*(R)\rangle$ converges to 0 as $\lambda$ tends to infinity. 
We apply Lemma \ref{lemma_star_exact} to our situation. We also use the fact that $T$  is closed and hence vanishes on exact test forms.
We obtain for every fixed $t\in\C^*$ that $\langle T, (A_{t\lambda})^*(R)\rangle-\langle T, (A_\lambda)^*(R)\rangle$ converges to 0. It follows that $\langle T_{t\lambda},R\rangle - \langle T_\lambda,R\rangle$ tends to 0. Therefore, tangent currents to $T$ along $V$ are invariant under the action of $(A_t)_*$, i.e. they are $V$-conic currents.

We prove that the tangent currents to $T$ have the same cohomology class. Let $S$ be such a current. Fix also a small neighbourhood of $V$ in $E$. It is not difficult to see that for $\lambda$ large enough $T_\lambda$ restricted to this neighbourhood is a closed current whose cohomology class does not depend on $\lambda$. It follows that  $\{S\}\smallsmile \{V\}$ does not depend on the choice of $S$.
Since $S$ is  $V$-conic and $V$ is a hypersurface, we deduce that the class of $S$ is $H^*_c(\overline E,\C)$ does not depend on the choice of $S$, see also Proposition \ref{prop_conic_char}. 
\hfill $\square$
\medskip

Let $S$ be a tangent current to $T$ along $V$. Denote by $\kappa^V(T)$
the class of $S$ in $H^{2p}_c(\overline E,\C)$. We know that it does not depend on the choice of $S$. 

\begin{definition}\rm \label{def_tangent_class}
We say that $\kappa^V(T)$ is the {\it total tangent class} of $T$ along $V$.
The h-dimension of $\kappa^V(T)$ is {\it the tangential h-dimension of $T$ along $V$}.  The {\it set of tangent directions} of $T$ along $V$ is the union of $\supp(S)$ for $S$ varying on the set of all tangent currents. Its projection to $V$ is the {\it tangent locus} of $T$ along $V$. 
\end{definition}

If $-h_{\overline E}$ denotes the tautological class of $\overline E$, as in Section \ref{section_fibration}, 
we can write in a unique way 
$$\kappa^V(T)=\sum_{j=\max(0,l-p)}^{\min(l,k-p)} \pi_0^*(\kappa_j^V(T))\smallsmile h_{\overline E}^{j-l+p},$$
where $\kappa_j^V(T)$ is a class in $H^{2l-2j}_c(V,\C)$.  By convention, $\kappa_j^V(T)$ is 0 when $j$ does not satisfy the inequalities $\max(0,l-p)\leq j\leq \min(l,k-p)$. With notations as in Section \ref{section_fibration}, we have $\kappa_j^V(T)=\kappa_j(S)$ if $S$ is a tangent current to $T$ along $V$.

\begin{remarks} \rm \label{rk_reduce_dim}
Let $\theta$ be a smooth positive closed $(q,q)$-form on $X$ with $q\leq k-p$. Let $S$ be the tangent current to $T$ associated with a sequence $(\lambda_n)$. Then $S\wedge\pi_0^*(\theta_{|V})$ is the tangent current to $T\wedge \theta$ associated with the same sequence. We also have $\kappa^V(T\wedge\theta)=\kappa^V(T)\smallsmile \pi_0^*\{\theta_{|V}\}$ and $\kappa_j^V(T\wedge\theta) = \kappa_{j+q}^V(T)\smallsmile\{\theta_{|V}\}$. 
\end{remarks}

We consider now a case which is very useful in computing tangent classes. 

\begin{lemma} \label{lemma_hyp_case}
Let $X,V$ and the $(p,p)$-current $T$ be as above.
Assume  that $p\leq l$ and that the tangential h-dimension of $T$ along $V$ is minimal, i.e. equal to $l-p$. Then 
$\kappa_{l-p}^V(T)=\{T\}_{|V}$ and $\kappa^V(T)=\pi_0^*(\{T\}_{|V})$. In particular, 
when $V$ is a hypersurface of $X$ and $T$ has no mass on $V$, 
the above identities hold and we have moreover $\kappa^V(T)_{|H_\infty}=\{T\}_{|V}$.
\end{lemma}
\proof
When $V$ is a hypersurface,  we have $\P(E)=V$ and $V$-conic currents without mass on $V$ are pull-back by $\pi_0$ of currents on $V$. Moreover, $\pi_0$ defines an isomorphism between $H_\infty$ and $V$. Therefore, the second assertion is a direct consequence of the first one.
We prove the first assertion using the notation introduced above. For simplicity, assume that $\tau$ is smooth.

By Lemma \ref{lemma_bi_lipschitz}, we have 
$\{T_\lambda\}_{|V}=\{T\}_{|V}$ for every $\lambda$. So if $S$ is a tangent current to $T$ along $V$, the class  $\{S\}_{|V}$  is equal to $\{T\}_{|V}$. Lemma \ref{lemma_h_dim_min}  implies that $S$ is the pull-back by $\pi_0$ of  the shadow $S^h$ of $S$ on $V$. Therefore, we have $\kappa^V(T)=\{S\}=\{\pi_0^*(S^h)\}$ and  $S^h$ belongs to the class $\{S\}_{|V}=\{T\}_{|V}$. 
On the other hand, by Lemma \ref{lemma_current_class}, $S^h$ belongs to
$\kappa_{l-p}(S)=\kappa_{l-p}^V(T)$. The lemma follows.
\endproof

The following result shows the upper semi-continuity for the maximal h-dimensional part of the tangent class when the current $T$ varies. 

\begin{theorem} \label{thm_tangent_usc}
Let $X,V$ and $W_1$ be as in Theorem \ref{th_tangent_conic}.
Let $T_n$ and $T$ be positive closed $(p,p)$-currents on $X$ with support in $W_1$ such that $T_n\to T$. Let $s$ be the tangential h-dimension of $T$ along $V$. Then
\begin{enumerate}
\item  If $r$ is an integer strictly larger than $s$, then $\kappa_r^V(T_n)$ converges to $0$.
\item If $\kappa_s$ is a limit class of the sequence $\kappa_s^V(T_n)$, then the classes $\kappa_s$  and 
$\kappa_s^V(T)-\kappa_s$ are pseudo-effective. 
\end{enumerate}
\end{theorem}
\proof
If $T$ has positive mass on $V$, then the tangential h-dimension of $T$ along $V$ is maximal, i.e. equal to $k-p$. 
The theorem is clear. Assume now that $T$ has no mass on $V$. We deduce that the mass of $T_n$ on $V$ tends to 0. So removing from $T_n$ its restriction to $V$ permits to assume that $T_n$ has no mass on $V$ for every $n$. 
 
Denote by $\widehat T$ and $\widehat T_n$ the strict transforms of $T$ and $T_n$ with respect to the blow-up $\sigma:\widehat X\to X$ along $V$. 
Recall that we identify the hypersurface at infinity $H_\infty$ of $\overline E$ with $\sigma_E^{-1}(H_\infty)$ and with $\P(E)$. So the restriction of a class $\kappa$ to $H_\infty$ or to $\sigma_E^{-1}(H_\infty)$ is denoted by $\kappa_{|\P(E)}$. 
By Lemma \ref{lemma_tangent_blowup} and the last assertion of Lemma \ref{lemma_hyp_case}, we have 
$$\kappa^V(T)_{|\P(E)}=\kappa^{\widehat V}(\widehat T)_{|\P(E)}=\{\widehat T\}_{|\widehat V}$$ 
and a similar property for $T_n$. 
Extracting a subsequence we can assume that 
$\widehat T_n$ converges to a current $\widehat T'$. Write $\widehat T' =\widehat T+\widehat R$ where   $\widehat R$ is the restriction of $\widehat T'$ to $\widehat V$. If $\{\widehat R\}'$ denotes the class of $\widehat R$ in $H^{2p}_{\widehat V}(\widehat X,\C)$, we have
$$\lim_{n\to\infty} \kappa^{V}(T_n)_{|\P(E)}-\kappa^V(T)_{|\P(E)}=\{\widehat R\}'_{|\widehat V}.$$
We show that the h-dimension of $\widehat R$ is at most equal to $s$. 

Assume that the h-dimension of $\widehat R$ is strictly larger than $s$. By Remarks \ref{rk_reduce_dim}, we can multiply $T$ and $T_n$ by a strictly positive closed form in order to assume that $\kappa^V(T)=0$ and that $\widehat R$ is a vertical current. 
We then have $\lim \kappa^{V}(T_n)_{|\P(E)}=\{\widehat R\}'_{|\widehat V}$. 
On one hand the above limit is a pseudo-effective class.  On the other hand,
arguing as in Lemma \ref{lemma_negative_intersection}, the class $\{\widehat R\}'_{|\widehat V}$ can be represented by a strictly negative constant times a linear subspace on a fiber of $\pi$. This is a contradiction.

So the $h$-dimension of $\widehat R$ is at most equal to $s$. We deduce from the above computation on $\lim \kappa^{V}(T_n)_{|\P(E)}$ that the h-dimension of this limit is at most equal to $s$. This gives us the first part of the theorem, see Lemma \ref{lemma_compare_infty}. Since the last limit is a pseudo-effective class, we also deduce that
$\lim\kappa_s^{V}(T_n)$ is pseudo-effective. This implies that the class $\kappa_s$ in the second part of the theorem is pseudo-effective. Finally, by Lemma \ref{lemma_negative_intersection}, the class $-\kappa_s(\{\widehat R\}'_{|\widehat V})$ is pseudo-effective. This and the above computation imply that $\kappa_s^V(T)-\kappa_s$ is pseudo-effective and complete the proof of the theorem.
\endproof

Note that when $T$ has no mass on $V$ its total tangent class $\kappa^V(T)$ along $V$ is determined by its restriction to the hypersurface at infinity $H_\infty\simeq \P(E)$. As above, we denote this class by $\kappa^V(T)_{|\P(E)}$. 
We identify both $\pi_0:H_\infty\to V$ and  $\sigma:\widehat V\to V$ with $\pi:\P(E)\to V$. 
The following proposition gives us a way to compute the tangent class of $T$ along $V$. It is similar to Siu's point of view on the Lelong number at a point using the blow-up at this point, see Siu \cite{Siu}. 

\begin{proposition} Let $X,V$ and $T$ be as in Theorem \ref{thm_tangent_usc}. Let $\widehat T$ be the strict transform of $T$ with respect to the blow-up $\sigma:\widehat X\to X$ of $X$ along $V$. Denote by $-h_{\P(E)}$ the tautological class of $\pi:\P(E)\to V$ as above. Let $e(T)$ be the class in $H^{2p-2}_c(\P(E),\C)$ defined in Lemma \ref{lemma_exceptional_class}.
Assume that $T$ has no mass on $V$. 
Then 
$$\kappa^V(T)_{|\P(E)}=e(T)\smallsmile h_{\P(E)} +\pi^*(\{T\}_{|V}).$$
\end{proposition}
\proof
For simplicity, we  identify $\sigma_E^{-1}(H_\infty)$ with $H_\infty$ and with $\P(E)$.
We then have $\kappa^{ V}( T)_{|\P(E)}=\kappa^{\widehat V}(\widehat T)_{|\P(E)}$. 
By Lemma \ref{lemma_hyp_case} applied to $\widehat V$, we have 
$\kappa^{\widehat V}(\widehat T)_{|\P(E)}=\{\widehat T\}_{|\widehat V}.$ 
Recall that $\{\widehat T\}$ is equal to the difference between $\sigma^*\{T\}$ and the canonical image $\widetilde e(T)$ of $e(T)$ in $H^{2p}_{\widehat V}(\widehat X,\C)$. 
We also have 
$(\sigma^*\{T\})_{|\widehat V}=\pi^*(\{T\}_{|V})$ (this can be seen using a smooth form in $\{T\}$). 
Moreover, 
$$\widetilde e(T)_{|\widehat V}=\big(\pi_0^*(e(T))\smallsmile [\widehat V]\big)_{|\widehat V}=e(T)\smallsmile \{\widehat V\}_{|\widehat V}=-e(T)\smallsmile h_{\P(E)}.$$
This  implies the proposition. Note that $\{T\}_{|V}=0$ when $p>\dim V$. 
\endproof

The following result will be used to bound tangent classes and to show that some tangent classes vanish.

\begin{proposition} \label{prop_cone_restriction}
Let $X,V$ and $T$ be as above. 
Let $V'$ be a submanifold of $V$. Let $s$ denote the tangential h-dimension of $T$ along $V$. 
Then the tangential h-dimension of $T$ along $V'$ is at most equal to $s$. Moreover, if $S$ is a tangent current to $T$ along $V$, we have $\kappa^{V'}_s(T)\leq \kappa^{V'}_s(S)$. The inequality still holds if we replace $s$ with the tangential h-dimension of $S$ along $V'$. 
\end{proposition}

If $T$ has support in $V$, then $S=T$ and the proposition is clear. So we can assume that $T$ has no mass on $V$. In particular, we have $s<k-p$.

Let $\tau$ be the smooth admissible map given in Lemma \ref{lemma_admissible}. 
Let $\sigma':\widehat X'\to X$ be the blow-up of $X$ along $V'$ and $\sigma_{E'}:\widehat E'\to E'$ the blow-up along $V'$ of the normal vector bundle $E'$ to $V'$ in $X$. 
Let $\widehat T'$ be the strict transform of $T$ by $\sigma':\widehat X'\to X$.
Define $\widehat V':=\sigma'^{-1}(V')$ and we identify it with $\sigma_{E'}^{-1}(V')$ and also with $\P(E')$. 
Observe that in general $\tau$ is not admissible with respect to $V'$. 
We need the following lemma.

\begin{lemma} The map $\tau$ lifts to a bi-Lipschitz map $\widehat\tau'$ from a neighbourhood of $\widehat V'$ in $\widehat X'$ to a neighbourhood of $\widehat V'$ in $\widehat E'$ which is smooth outside $\widehat V'$ and preserves the hypersurface $\widehat V'$. Moreover, if $\widetilde T':=\widehat\tau'_*(\widehat T')$,
we have $\{\widehat T'\}_{|\widehat V'}=\{\widetilde T'\}_{|\widehat V'}$.
\end{lemma}
\proof
We have $\widehat \tau'=\sigma_{E'}^{-1}\circ\tau\circ\sigma'$ outside $\widehat V'$. We first show that this map extends to a bi-Lipschitz map. The map $\tau$ is described locally as in Definition \ref{def_admissible_loc} and in the proof of Lemma \ref{lemma_admissible} where all functions involved are smooth.  In order to simplify the notation, we will not use exactly the same coordinate system of $\D^k$ as above. 

Let $(y^1,y^2)$ denote a linear coordinates system on $\D^k$ where $y^1=(y_1,\ldots,y_{l'}):=x^1$ but $y^2$ is obtained from $(x^2,x^3)$ by an index permutation. The aim is to consider that the components of $y^2$ play an equivalent role. We can write in these coordinates
$$\tau(y)=\big(y^1+y^2b(y),y^2c(y)\big)+O(\|y^2\|^2) \quad \mbox{as} \quad y^2\to 0,$$
where the functions involved in $b,c$ and $O(\|y^2\|^2)$ are smooth and the determinant of the matrix $c(y)$ is equal to 1.
In these coordinates we identify $\sigma'$ with $\sigma_{E'}$.  

We cover $\sigma'^{-1}(\D^k)$ with a finite number of equivalent charts and as above we will only work in one of them. The considered chart is denoted by $\widehat D'$ endowed with coordinates
$w=(w^1,w^\#,w_k)$ with $w^1:=(w_1,\ldots,w_{l'})$, $w^\#:=(w_{l'+1},\ldots,w_{k-1})$, $|w_j|<2$ such that
$$\sigma'(w)=(w^1,w_kw^\#,w_k) \quad \mbox{and} \quad \sigma_{E'}^{-1}(y)=\sigma'^{-1}(y)=(y^1,y_k^{-1}y^\#,y_k).$$
We deduce that
$$\widehat\tau'(w)=\big(w^1+w_k\widetilde b(w)+O(|w_k|^2),\widetilde c^\#(w)+w_k^{-1}O(|w_k|^2),w_k\widetilde c_k(w)+O(|w_k|^2)\big),$$
where the functions involved in $\widetilde b$, $\widetilde c^\#$, $\widetilde c_k$ and $O(|w_k|^2)$ are smooth. 
The inverse of $\widehat\tau'$ satisfies a similar property.
We see that $\widehat\tau'$ extends to a bi-Lipschitz map which is not identity on $\widehat V'$ in general.
The hypersurface $\widehat V'$ is given by $w_k=0$. So it is invariant under $\widehat\tau'$. 

It remains to prove the last identity in the lemma. By Lemma \ref{lemma_bi_lipschitz}, we only have to check that the restriction $\widetilde\tau$ of $\widehat\tau'$ to $\widehat V'$ acts trivially on $H_c^*(\widehat V',\C)$. In local coordinates as above, we have  $\widetilde\tau(w^1,w^\#)= (w^1,\widetilde c^\#(w^1,w^\#,0))$. So it is induced by the differential of $\tau$ which is $\C$-linear at each point of $V$. We deduce that $\widetilde\tau$ is induced by a smooth self-map of the tautological line bundle $O_{\widehat V'}(-1)$ of $\widehat V'$ which sends $\C$-linearly fibers to fibers. It follows that $\widetilde\tau$ preserves the tautological class of $\widehat V'$. On the other hand, it preserves the fibers over $V'$. Hence, Leray's theorem implies that $\widetilde\tau$ acts trivially on $H_c^*(\widehat V',\C)$. This completes the proof of the lemma.
\endproof

\noindent
{\bf End of the proof of Proposition \ref{prop_cone_restriction}.} 
Recall that $T$ has no mass on $V$ and $s<k-p$. 
By Lemma \ref{lemma_hyp_case} applied to $\widehat V'$, we have $\kappa^{V'}(T)_{|\P(E')}=\{\widehat T'\}_{|\widehat V'}$. It follows from the last lemma that $\kappa^{V'}(T)_{|\P(E')}=\{\widetilde T'\}_{|\widehat V'}$. The map $A_\lambda$ can be lifted to a holomorphic map $\widehat A_\lambda:\widehat E'\to\widehat E'$. Since this map depends continuously on $\lambda$, it acts trivially on cohomology with integer coefficients. Therefore, it acts trivially on de Rham cohomology. Thus, $\kappa^{V'}(T)_{|\P(E')}=\{\widehat T_\lambda'\}_{|\widehat V'}$ where $\widehat T_\lambda':=(\widehat A_\lambda)_*(\widetilde T')$.  
Define $T_\lambda:=(A_\lambda)_*\tau_*(T)$ as above. We have $(\sigma_{E'})_*(\widehat T_\lambda')=T_\lambda$. 
Let $(\lambda_n)$ be a sequence such that $T_{\lambda_n}$ converges to $S$. 
Then, $\widehat T_{\lambda_n}'$ converges outside $\widehat V'$ to the strict transform $\widehat S'$ of $S$ by $\sigma_{E'}$. We show that any limit current of $\widehat T_{\lambda_n}'$ is equal to $\widehat S'$ plus a positive closed current supported by $\widehat V'$.

Let $R$ be a smooth $(2k-2p)$-form with compact support in $\widehat X'$. We show that the family of $\langle \widehat T_\lambda',R\rangle$ is bounded for $\lambda$ large enough. Using a partition of unity, we reduce the problem to the case where $R$ is supported by $\sigma'^{-1}(\D^k)$ as in Lemma
\ref{lemma_star_admissible_bis}. Since the considered currents have no mass on $\widehat V'$ and $V'$, we have
$$\langle \widehat T_\lambda',R\rangle = \langle T_\lambda,\sigma'_*(R)\rangle_{\D^k\setminus V'}=\big\langle T,\tau^*(A_\lambda)^*\sigma'_*(R)\big\rangle_{\D^k\setminus V'}.$$
Lemmas \ref{lemma_mass_star} and \ref{lemma_star_admissible_bis}  imply that the family $\langle \widehat T_\lambda',R\rangle$ is bounded. It follows that the family of currents $\widehat T_\lambda'$ is relatively compact.

By Lemma \ref{lemma_star_admissible}, if the component of bidegree $(k-p,k-p)$ of  $R$ vanishes, the above integral converges to 0. Therefore, the limit currents of $\widehat T_\lambda'$ are of bidegree $(p,p)$. The same proposition shows that if $R$ is a weakly positive $(k-p,k-p)$-form, then the limit values of 
$\langle \widehat T_\lambda',R\rangle$ is positive. We conclude that the limit currents of $\widehat T_{\lambda_n}'$ are positive closed $(p,p)$-currents. Recall that these currents are equal to $\widehat S'$ outside $\widehat V'$. Let $\widehat S'+\widehat S''$ be such a limit current with $\widehat S''$ positive closed supported by $\widehat V'$. 
Denote by $\{\widehat S''\}'$ the class of $\widehat S''$ in $H^{2p}_{\widehat V'}(\widehat X',\C)$. 

We deduce from the above discussion that 
$$\kappa^{V'}(T)_{|\P(E')}=\{\widehat S'\}_{|\widehat V'}+\{\widehat S''\}'_{|\widehat V'}=\kappa^{V'}(S)_{|\P(E')}+\{\widehat S''\}'_{|\widehat V'}.$$
Let $r$ denote the h-dimension of $T$ along $V'$. If $r$ is strictly larger than $s$, replacing $T$ by $T\wedge\omega^r$ gives us identities similar to the last ones with $S=0$; this contradicts Lemma \ref{lemma_negative_intersection} applied to $\widehat S''$ and the fact that $\kappa_r^{V'}(T)_{|\P(E')}$ is a non-zero pseudo-effective class. So we have $r\leq s$.  
The same property applied to $S$ in $\overline E$ says that 
the tangential h-dimension $s'$ of $S$ along $V'$ is at most equal to $s$. Therefore, using the above identities again, we conclude that 
the h-dimension $s''$ of $\widehat S''$ with respect to $V'$ is at most equal to $s$ and we also have 
$r\leq \max(s',s'')$. 

If $s'<s''$, a similar argument as above for $T\wedge \omega^{s''}$ gives a contradiction. Therefore, we have $s''\leq s'$ and $r\leq s'$. 
Then, using Lemmas \ref{lemma_compare_infty} and  \ref{lemma_negative_intersection} we obtain that
$$\kappa_{s'}^{V'}(T)=\kappa^{V'}_{s'}(S)+\kappa_{s'}(\{\widehat S''\}'_{|\widehat V'})\leq \kappa^{V'}_{s'}(S).$$ 
This is the last assertion in the proposition.

If $s=s'$, the first assertion in the proposition is clear. Otherwise, we have $r\leq s' <s$ and the assertion is also clear because 
$\kappa_{s}^{V'}(T)=\kappa^{V'}_s(S)=0$. This completes the proof of the proposition. 
\hfill $\square$

\section{Density and intersection of currents} \label{section_density}

Let $X$ be a K\"ahler manifold of dimension $k$ as above. 
In this section we will introduce a notion of density associated with any finite family of positive closed currents such that the intersection of their supports is compact. The last condition is satisfied when $X$ is already compact. 
We will study some basic properties of the density and compare it with the Lelong number. We also discuss a new notion of  intersection of currents and compare it with classical notions.

Let $T_j$ be a positive closed current of bidegree $(p_j,p_j)$ on $X$ with $1\leq j\leq m$. Assume that the intersection of their supports is compact. Define $\T:=T_1\otimes \cdots\otimes T_m$. This is a positive closed $(p,p)$-current on $X^m$ with $p:=p_1+\cdots+p_m$. Denote by $\Delta$ the diagonal of $X^m$, i.e. the set of points $(x,\ldots,x)$ with $x\in X$. It is canonically isomorphic to $X$. Then the intersection of $\supp(\T)$ with $\Delta$ is compact
and $\T$ has no mass on $\Delta$ except when the $T_j$ are measures which contain a same atom.

Denote by $\Tan(X)$, $\Tan(X^m)$ and $\Tan(\Delta)$ the tangent vector bundles of $X,X^m$ and $\Delta$ respectively.
Let $\E_m$ denote the normal bundle to $\Delta$ in $X^m$. 
The vectors which are tangent to the fibers of 
the natural projection $(x_1,\ldots, x_m) \mapsto (x_1,\ldots,x_{m-1})$ constitute a vector sub-bundle of $\Tan(X^m)$. Its restriction to $\Delta$ is a complement of $\Tan(\Delta)$ in $\Tan(X^m)$. We see that $\E_m$ is canonically isomorphic to $\Tan(X)\oplus\cdots\oplus \Tan(X)$ ($m-1$ times). 
So the rank of $\E_m$ is equal to $(m-1)k$.
Define 
$$\kappa(T_1,\ldots,T_m):=\kappa^\Delta(\T).$$ 
This is a pseudo-effective cohomology class in $H_c^{2p}(\overline\E_m,\C)$.
Define also 
$$\kappa_j(T_1,\ldots,T_m):=\kappa_j^\Delta(\T).$$ 
This is a cohomology class in $H^{2k-2j}_c(X,\C)$.

\begin{definition}\rm
The class $\kappa(T_1,\ldots,T_m)$ is called {\it the total density class}; 
the class $\kappa_j(T_1,\ldots,T_m)$ is {\it the density class of (complex) dimension $j$} and
the h-dimension of $\kappa(T_1,\ldots,T_m)$ is {\it the density h-dimension}  associated with $T_1,\ldots,T_m$. 
If $\S$ is a tangent current to $\T$ along $\Delta$, we say that $\S$ is {\it a density current associated with $T_1,\ldots,T_m$}. 
\end{definition}

Observe that any permutation of $(x_1,\ldots,x_m)$ induces holomorphic automorphisms of $\E_m$ and of $\overline \E_m$ which leave invariant the fibers. So
this action on $\overline \E_m$ preserves the tautological class $-h_{\overline\E_m}$. 
We then deduce from Leray's theorem that the action is in fact the identity on $H^*_c(\overline\E_m,\C)$.  Therefore,  
$\kappa$ and $\kappa_j$ are symmetric in $T_1,\ldots,T_m$.

\begin{example} \rm
If the currents $T_j$ have locally continuous potentials, we can show that the current $\T$ admits a unique tangent current along $\Delta$. This current vanishes when $p_1+\cdots+p_m>k$ and is equal to the pull-back of $T_1\wedge\ldots\wedge T_m$ otherwise, see also Proposition \ref{prop_compare_wedge} below.
\end{example}

\begin{lemma} \label{lemma_density_h_dim}
The density h-dimension associated with $T_1,\ldots,T_m$ is smaller or equal to the complex dimension $k-p_j$ of $T_j$ for $1\leq j\leq m$. 
\end{lemma}
\proof
Let $s$ denote the density h-dimension  associated with $T_1,\ldots,T_m$. The lemma is clear if 
the density class  vanishes. 
Suppose this is not the case. 
Then the class $\kappa_s(T_1,\ldots,T_m)$ is non-zero and pseudo-effective. We have
$$\kappa_s(T_1,\ldots,T_m)\smallsmile \{\omega^s\}\not=0.$$
The last class is also the shadow of $\kappa^\Delta(T_1\otimes\cdots\otimes T_m\wedge \Pi_j^*(\omega^s))$, where $\Pi_j$ is the projection from $X^m$ to the $j$-th factor.  We deduce that $T_1\otimes\cdots\otimes T_m\wedge \Pi_j^*(\omega^s)\not =0$.  
It follows that $T_j\wedge\omega^s\not=0$ and hence $s\leq k-p_j$. 
\endproof

The following lemma shows that the notion of density generalizes the notion of tangent currents.

\begin{lemma}
Let $X,V$ and $T$ be as in Section \ref{section_tangent}. Then $\kappa_j(T,[V])$ is equal to the canonical image of $\kappa^V_j(T)$ in $H^{2k-2j}_c(X,\C)$. 
\end{lemma}
\proof
Observe that 
the pull-back of $T$ by the canonical projection $\Pi:V\times X \to X$ can be identified to 
the current $[V]\otimes T$ in $X\times X$. The restriction of $\E_2$ to $\Delta_V:=(V\times X)\cap \Delta$ can be identified with the pull-back of the tangent vector bundle of $X$ by the restriction $\Pi_{|\Delta_V}$ of $\Pi$ to $\Delta_V$. We denote it by $F$. 
The tangent currents to $[V]\otimes T$ along $\Delta$ can be identified with tangent currents to $\Pi^*(T)$ along $\Delta_V$. 
The pull-back to $\Delta_V$ by $\Pi$ of the tangent vector bundle of $V$  is a sub-bundle of $F$ that we denote by $F'$. The quotient $F/F'$ can be identified with the normal vector bundle $E$ to $V$ in $X$ if we identify $\Delta_V$ with $V$. Denote by $\rho:F\to F/F'$ the canonical projection. 

We show that  the tangent currents to $\Pi^*(T)$ along $\Delta_V$ are equal to the pull-back by $\rho$ of the tangent currents to $T$ along $V$. For this purpose, we will use local coordinates as in Section \ref{section_tangent}. 
We identify a chart of $X$ with $\D^k=\D^l\times \D^{k-l}$ on which $V$ is equal to $\D^l\times\{0\}$. Consider the natural  coordinate system $(x',y',y'')$ on the chart $\D^l\times\D^k$ of $V\times X$ where $\Delta_V$ is given by $\{y'=x', y''=0\}$ and $\Pi(x',y',y'')=(y',y'')$.  

In order to identify $\Delta_V$ with $V$ we use the coordinate system $(x',z',y'')$ with $z':=y'-x'$. So $V$ is identified with $\Delta_V$ and given by $\{z'=0, y''=0\}$. The vector bundle $F$ is then identified to $\D^l\times\C^k$ and $F'$ is the intersection of $\D^l\times \C^k$ with the subspace $\{y''=0\}$. So it is equal to $\D^l\times\C^{l}\times \{0\}$. The vector bundles $F/F'$ and $E$ are identified to $\D^l\times \{0\}\times \C^{k-l}$. The map $\rho$ is just the canonical projection $(x',z',y'')\mapsto (x',0,y'')$.

The projection $\Pi$ is given in these coordinates by $\Pi(x',z',y'')=(x'+z',0,y'')$ where we identify the chart $\D^k\subset X$ with the polydisc $\D^l\times\{0\}\times\D^{k-l}$. 
We use the identity map as a (local) admissible map associated with $X$ and $V$. We also use the map $\tau(x',z',y''):=(x'+z',z',y'')$ for the pair $V\times X$ and $\Delta_V$. The multiplication with $\lambda$ on $E$ and $\E_2$ are identified with the map $A_\lambda(x',z',y''):=(x',\lambda z',\lambda y'')$.  
We see that $\Pi\circ (A_\lambda\circ\tau)^{-1}=(A_\lambda)^{-1}\circ \rho$. It follows that 
$(A_\lambda)_*\tau_*\Pi^*(T)=\rho^*(A_\lambda)_*(T)$. Thus, if $S$ is a tangent current to $T$ along $V$ then $\rho^*(S)$ is a tangent current to $\Pi^*(T)$ along $\Delta_V$. 

The map $\rho$ induces a meromorphic map $\widetilde\rho:\overline \E_2\to \overline E$. It is not difficult to show that if $-h_{\overline E}$ is the tautological class of $\overline E$ then $-\widetilde\rho^*(h_{\overline E})$ is the tautological class of $\overline E$.  Finally, the uniqueness of the decomposition in Leray's theorem, implies that $\kappa_j(\widetilde \rho^*(S))=\kappa_j(S)$. The lemma follows.

Note that the above construction gives an isomorphism between the set of tangent currents to $T$ along $V$ and the set of tangent currents to $[V]\otimes T$ along the diagonal $\Delta$. 
\endproof

Note that in the last lemma the canonical morphism from $H^{2l-2j}_c(V,\C)$ to $H^{2k-2j}_c(X,\C)$ is not injective in general. However, the lemma still holds if we replace $X$ by a small enough neighbourhood of $V$ and in that case the corresponding morphism is injective.

Let $T$ be a positive closed $(p,p)$-current on $X$. By Siu's theorem \cite{Siu}, the Lelong number $\nu(T,x)$ defines a function which is upper semi-continuous with respect to the Zariski topology on $X$. In particular, if $Y$ is an irreducible analytic set, then $\nu(T,\cdot)$ is constant on a dense Zariski open set of $Y$. We denote this constant by $\nu(T,Y)$. Moreover, also by Siu's theorem,
there is a finite or countable family of irreducible analytic sets $Y_j$ of dimension $k-p$ and constants $c_j> 0$ such that 
$$T=\sum c_j[Y_j]+T'$$
where $T'$ is a positive closed current such that for every $c>0$ the level set $\{\nu(T,\cdot)\geq c\}$ is an analytic set of dimension $\leq k-p-1$. The following results give the relation between density of currents and Lelong numbers.

\begin{lemma} \label{lemma_density_lelong}
Let $X$ and $T_j$ be as above. Assume that $T_1$ is a measure, i.e. $p_1=k$. Then the density h-dimension  associated with $T_1,\ldots,T_m$ is $0$ and we have
$$\kappa_0(T_1,\ldots,T_m)=\big\langle T_1,\nu(T_2,\cdot)\ldots \nu(T_m,\cdot)\big\rangle.$$
In particular, the function $a\mapsto \kappa_0(\delta_a,T_2,\ldots,T_m)$ is upper semi-continuous with respect to the Zariski topology on $X$.
\end{lemma}
\proof
The first assertion is a consequence of Lemma \ref{lemma_density_h_dim}. For the second assertion, since $\kappa_0$ is linear on each variable, we can disintegrate $T_1$ into Dirac masses and assume for simplicity that $T_1$ is the Dirac mass at a point $a$. In this case, we see that $\kappa_0(T_1,\ldots,T_m)$ is the Lelong number of $T_2\otimes \cdots\otimes T_m$ at the point $(a,\ldots,a)$. It is not difficult to see that this Lelong number is equal to 
$\nu(T_2,a)\ldots\nu(T_m,a)$. The lemma follows.
\endproof

\begin{proposition} \label{prop_density_lelong}
Let $X$ and $T_j$ be as above. Assume that $T_1$ is 
the current of integration on an irreducible analytic set $Y$ of dimension $k-p_1$. Then 
$$\kappa_{k-p_1}(T_1,\ldots,T_m)=\nu(T_2,Y)\ldots\nu(T_m,Y)\{Y\}.$$
\end{proposition}
\proof
Let $\S$ be a density current associated with $T_1,\ldots,T_m$. If $\kappa_{k-p_1}(T_1,\ldots,T_m)$ does not vanish, by Lemma \ref{lemma_density_h_dim}, the h-dimension of $\S$ is equal to $k-p_1$. So the last class contains the shadow of $\S$ which is a positive closed $(p_1,p_1)$-current supported by $Y$. Therefore, it is equal to a constant $c$ times $\{Y\}$.  Of course, this property holds also when the considered class vanishes. We compute now the constant $c$.

We have
$$\kappa_0(T_1\wedge\omega^{k-p_1},T_2,\ldots,T_m)=\kappa_{k-p_1}(T_1,\ldots,T_m)\smallsmile \{\omega^{k-p_1}\}
=c\{T_1\wedge\omega^{k-p_1}\}.$$
Since $T_1\wedge\omega^{k-p_1}$ is a positive measure, by Lemma \ref{lemma_density_lelong}, we have 
$$c\{T_1\wedge\omega^{k-p_1}\}=\big\langle  T_1\wedge\omega^{k-p_1}, \nu(T_2,\cdot)\ldots \nu(T_m,\cdot)\big\rangle.$$
Since $T_1=[Y]$, Siu's theorem mentioned above implies that the last integral is equal to $\nu(T_2,Y)\ldots \nu(T_m,Y) \{T_1\wedge\omega^{k-p_1}\}$. The proposition follows.
\endproof

\begin{proposition}
Let $X$ and $T_j$ be as above. 
Assume that the set $\Ec$ of points $x$ such that $\{\nu(T_2,x)>0\}$ contains no analytic set of dimension $k-p_1$. Then the density h-dimension associated with $T_1,\ldots,T_m$ is strictly smaller than $k-p_1$. 
\end{proposition}
\proof
By Lemma \ref{lemma_density_h_dim}, this dimension is at most equal to $k-p_1$. So it is enough to prove that $\kappa_{k-p_1}(T_1,\ldots,T_m)\smallsmile \{\omega^{k-p_1}\}=0$. Arguing as in the end of Proposition \ref{prop_density_lelong}, we have 
$$\kappa_{k-p_1}(T_1,\ldots,T_m)\smallsmile \{\omega^{k-p_1}\}=\big\langle T_1\wedge\omega^{k-p_1}, \nu(T_2,\cdot)\ldots \nu(T_m,\cdot) \big\rangle.$$
By hypothesis, $\Ec$ is a finite or countable union of analytic sets of dimension less that $k-p_1$. Therefore, the measure $T_1\wedge\omega^{k-p_1}$ has no mass on 
$\Ec=\{\nu(T_2,\cdot)\not =0\}$ and the last integral vanishes. The proposition follows.
\endproof

The following result is a direct consequence of Theorem \ref{thm_tangent_usc}.

\begin{corollary} 
Let $U_1,\ldots, U_m$ be open subsets with relatively compact intersection in a K\"ahler manifold $X$. Let $T_{j,n}$ and $T_j$ be positive closed $(p_j,p_j)$-currents with support in $U_j$ such that $T_{j,n}\to T_j$ as $n\to\infty$. Let $s$ denote the density h-dimension associated with $T_1,\ldots,T_m$. Then $\kappa_j(T_{1,n},\ldots,T_{m,n})\to 0$ for $j>s$. Moreover, any limit class of the sequence $\kappa_s(T_{1,n},\ldots,T_{m,n})$ is pseudo-effective and smaller or equal to $\kappa_s(T_1,\ldots,T_m)$. 
\end{corollary}

In what follows, we will introduce a new definition for the intersection of positive closed currents. We will give some basic properties needed in our dynamical application. We believe that the theory has an independent interest and has to be developed. Let $X$ and $T_i$ be as above. We assume that $p_1+\cdots+p_m\leq k$ which is a necessary condition to give a meaning to the intersection of the $T_j$.

\begin{definition}\rm \label{def_wedge}
Assume that the density h-dimension associated with $T_1,\ldots,T_m$ is minimal, i.e. equal to $k-p_1-\cdots-p_m$. Assume that there is a unique density current  $\S$ associated with $T_1,\ldots,T_m$. We define $T_1\curlywedge\ldots\curlywedge T_m$ as the shadow of $\S$ with respect to the fibration $\pi:\overline\E_m\to \Delta$.
\end{definition}

Observe that in this case $\T:=T_1\otimes \cdots \otimes T_m$ admits a unique tangent current $\S$ with respect to $\Delta$ and by Lemma \ref{lemma_h_dim_min}, it is equal to the pull-back of the current $T_1\curlywedge\ldots\curlywedge T_m$.  We deduce that
$$\{T_1\curlywedge\ldots\curlywedge T_m\}=\lim_{\lambda\to\infty} \{\T_\lambda\}_{|\Delta} = \{\T\}_{|\Delta},$$
where $\T_\lambda$ is defined as in Section \ref{section_tangent} for $\T,\Delta$ instead of $T,V$. It follows that
$$\{T_1\curlywedge\ldots\curlywedge T_m\}=\{T_1\}\smallsmile\cdots\smallsmile \{T_m\}.$$ 

The permutations of factors in $X^m$ induce bi-holomorphic self-maps on $\overline\E_m$ which preserve the fibers of $\pi_0:\overline \E_m\to X$. Since $S=\pi_0^*(T_1\curlywedge\ldots\curlywedge T_m)$ these bi-holomorphic maps also preserve $\S$. Hence,  the wedge-product $T_1\curlywedge\ldots\curlywedge T_m$ is symmetric with respect to $T_1,\ldots,T_m$. 

Let $T$ be a positive closed $(p,p)$-current and $T'$ a positive closed $(1,1)$-current with bounded local potentials.
We can define the wedge-product $T'\wedge T$ by
$T'\wedge T:=\ddc(uT)$ when $u$ is a local potential of $T'$ which is a bounded psh function, see e.g. 
\cite{Demailly3,FornaessSibony}.
The definition does not depend on the choice of $u$ and therefore extends to the global setting. This wedge-product gives a positive closed $(p+1,p+1)$-current. 

When $T_1,\ldots,T_{m-1}$ are of bidegree $(1,1)$ and have locally bounded potentials,  we can define 
$T_1\wedge\ldots\wedge T_m$ by induction. 
The following result compares this wedge-product with the definition given above. The proposition holds for currents with bounded local potentials but the proof in that case is more technical.

\begin{proposition} \label{prop_compare_wedge}
Let $X$ and $T_j$ be as above. Assume that $T_1,\ldots,T_{m-1}$ are of bidegree $(1,1)$ with continuous local potentials. Then we have
$$T_1\wedge \ldots\wedge T_m= T_1\curlywedge\ldots\curlywedge T_m.$$
\end{proposition}
\proof
Let $U$ be a small open set in $X$ that we identify with the unit polydisc in $\C^k$ using local coordinates 
$x=(x_1,\ldots,x_k)$. Then $U^m$ is identified with the unit polydisc in $(\C^k)^m$ and we denote the canonical coordinate system by $(x^1,\ldots,x^m)$. 
We will use other coordinates $(y^1,\ldots,y^m)$ on $U^m$ given by $y^m:=x^m$ and $y^j:=x^j-x^m$ for $1\leq j\leq m-1$. In these coordinates, $\Delta\cap U^m$ is given by the equations $y^1=\cdots=y^{m-1}=0$. The vector bundle $\E_m$ is identified over $\Delta\cap U^m$ with $(\C^k)^{m-1}\times U$ where the zero section, i.e. $\Delta\cap U^m$, is identified with $\{0\}\times U$. 

In these coordinates, the application $\tau:=\id$ is admissible.  Write $T_j=\ddc u_j$ with $u_j$ psh  on $U$. 
Define a psh function $\widetilde u_j$ on $U^m$ by $\widetilde u_j(y^1,\ldots,y^m):=u_j(y^j+y^m)$. The current $\widetilde T_j:=\ddc \widetilde u_j$ is the pull-back of $T_j$ to $U^m$ by the projection from $U^m$ onto the $j$-th factor. Using convolutions, we can approximate $u_j$ uniformly by smooth psh functions. We see that 
$\T:=T_1\otimes\cdots\otimes T_m$ is equal to $\widetilde T_1\wedge \ldots\wedge \widetilde T_m$ where $\widetilde T_m$ is the pull-back of $T_m$ by the canonical projection  $\Pi:(\C^k)^{m-1}\times U\to U$.

Let $A_\lambda$ denote the multiplication by $\lambda$ along the factor $(\C^k)^{m-1}$ in $(\C^k)^{m-1}\times U$. We have to show that $(A_\lambda)_*(\T)$ converges to $\Pi^*(T_1\wedge\ldots\wedge T_m)$ as $\lambda\to\infty$.
We have 
$$(A_\lambda)_*(\T)=\ddc (\widetilde u_1\circ A_\lambda^{-1})\wedge \ldots \wedge \ddc (\widetilde u_{m-1}\circ A_\lambda^{-1})\wedge \widetilde T_m.$$
Since $\widetilde u_j\circ (A_\lambda)^{-1}$ converges locally uniformly to $u_j\circ\Pi$, the last wedge-product converges to $\Pi^*(T_1\wedge \ldots \wedge T_m)$. This completes the proof of the proposition.

Note that the same proof also works for currents of higher bidegree which admit local continuous potentials. 
\endproof

We will also need the following lemma in the dynamical setting.

\begin{lemma}
Let $X,U_j,T_j$ be as above. Let $1<l <m$ be an integer. Assume that $U:=U_{l+1}\cap\ldots\cap U_m$ is relatively compact in $X$ and that 
$T_j$ is of bidegree $(1,1)$ and has local continuous potentials in a neighbourhood of $\overline U$ for $1\leq j\leq l$. Assume also that 
$T_{l+1}\curlywedge \ldots \curlywedge T_m$ exists. Then $T_1\curlywedge\ldots\curlywedge T_m$ exists and we have
$$T_1\curlywedge\ldots\curlywedge T_m=T_1\wedge \ldots\wedge T_l\wedge (T_{l+1}\curlywedge \ldots \curlywedge T_m).$$
\end{lemma}
\proof
Using the similar notation as in Proposition \ref{prop_compare_wedge}, we have 
$$(A_\lambda)_*(\T)=\ddc (\widetilde u_1\circ A_\lambda^{-1})\wedge \ldots \wedge \ddc (\widetilde u_l\circ A_\lambda^{-1})\wedge (A_\lambda)_*(T_{l+1}\otimes\cdots\otimes T_m).$$
Observe that $\widetilde u_j\circ A_\lambda^{-1}$ converges locally uniformly to $u_j\circ\Pi$ and by hypotheses $(A_\lambda)_*(T_{l+1}\otimes\cdots\otimes T_m)$ converges to $\Pi^*(T_{l+1}\curlywedge\ldots\curlywedge T_m)$.
We deduce that the right hand side of the last identity converges to 
$$\Pi^*\big(T_1\wedge \ldots\wedge T_l\wedge (T_{l+1}\curlywedge \ldots \curlywedge T_m)\big).$$
The lemma follows.
\endproof

\begin{remark} \rm
The hypothesis on the supports of the currents $T_j$ can be refined. We can extend the notion of tangent current along a manifold $V$ to currents which satisfy suitable regularity near the points of $V\setminus K$ for some compact subset $K$ of $V$. It would be also useful to compare the above notion with the notion introduced in \cite{DinhSibony10, DinhSibony11}. We believe that the new notion extends the previous one and is valid in a more general setting. 
\end{remark}

\begin{remark} \rm
A subset $\Ec$ of $\R^+$ is said to be {\it of density zero} if the Lebesgue measure of $\Ec\cap [0,n]$ is equal to $o(n)$ when $n\to\infty$. Let $T_\lambda$ be as in Section \ref{section_tangent}. Assume that for some  set $\Ec$ of zero  density the limit of $T_\lambda$ exists when $\lambda\to\infty$ and $|\lambda|\not\in\Ec$. We then say that the limit is {\it the essential tangent current to $T$ along $V$}. We can consider a notion of intersection of currents $T_1,\ldots,T_m$ by assuming the existence of the essential tangent current to $T_1\otimes\cdots\otimes T_m$ along $\Delta$. We can also use another measure instead of the Lebesgue measure on $\R^+$
 or take some average before considering the limit.  
\end{remark}

The following example illustrates an advantage of the above notion of intersection. 

\begin{example} \rm
Let $\pi_1$ denote the canonical projection from $\P^1\times\P^1$ onto the first factor. Let $a$ be a point in $\P^1$ and $\nu$ a positive measure on $\P^1$ having no mass at $a$. Consider two positive closed currents $T_1:=\pi_1^*(\delta_a)$ and $T_2:=\pi_1^*(\nu)$. It is not difficult to check that $T_1\curlywedge T_2=0$. If the local potentials of $\nu$ are equal to $-\infty$ at $a$, then the local potentials of $T_2$ are equal to $-\infty$ on the support of $T_1$. In this case, we cannot define $T_1\wedge T_2$ in the classical sense. 
\end{example}

T.-C. Dinh, UPMC Univ Paris 06, UMR 7586, Institut de
Math{\'e}matiques de Jussieu, F-75005 Paris, France. {\tt
  dinh@math.jussieu.fr}, {\tt http://www.math.jussieu.fr/$\sim$dinh}

\

\noindent
N. Sibony,
Universit{\'e} Paris-Sud, Math{\'e}matique - B{\^a}timent 425, 91405
Orsay, France. {\tt nessim.sibony@math.u-psud.fr}

\end{document}